\documentclass[12pt]{article}
\usepackage{amsfonts,mathrsfs,amssymb,amsthm,mathptm}
\usepackage{amsmath,amscd}
\usepackage{mathptm,pslatex}
\usepackage{multirow}
\usepackage{color}
\usepackage[dvips]{graphicx}
\usepackage[all]{xy}
\usepackage{graphicx}
\usepackage{fancyhdr}
\usepackage{url}
\usepackage{exscale}
\usepackage{relsize}
\usepackage{bbm}

\usepackage{indentfirst}
\numberwithin{equation}{section}

\DeclareFontFamily{U}{cal}{}
\DeclareFontShape{U}{cal}{m}{n}{<->cmsy10}{}

\DeclareSymbolFont{rcal}{U}{cal}{m}{n}
\DeclareSymbolFontAlphabet{\mathcal}{rcal}

\newtheorem{Def}{Definition}[section]
\newtheorem{Bsp}[Def]{Example}

\newtheorem{Theo}[Def]{Theorem}
\newtheorem{Lem}[Def]{Lemma}
\newtheorem{Koro}[Def]{Corollary}
\theoremstyle{definition}
\newtheorem{Rem}[Def]{Remark}

\newcommand{\bsm}{\begin{smallmatrix}}
\newcommand{\esm}{\end{smallmatrix}}

\newcommand{\add}{{\rm add}}
\newcommand{\gd}{{\rm gl.dim }}

\newcommand{\End}{{\rm End}}

\def\ed{\mathop{\rm ext.dim}\nolimits}
\def\wrd{\mathop{\rm w.resol.dim}\nolimits}

\newcommand{\DD}{{\rm D}}

\newcommand{\rad}{{\rm rad}}

\newcommand{\pd}{{\rm pd}}
\newcommand{\id}{{\rm id}}

\newcommand{\DTr}{{\rm DTr}}

\newcommand{\CC}{{\mathcal C}}
\def\E{\mathcal{E}}
\newcommand{\F}{{\mathcal F}}
\def\H{\mathcal{H}}

\newcommand{\T}{{\mathcal T}}
\newcommand{\U}{{\mathcal U}}
\newcommand{\V}{{\mathcal V}}

\newcommand{\cpx}[1]{#1^{\bullet}}
\newcommand{\D}[1]{{\mathscr D}(#1)}

\newcommand{\Db}[1]{{\mathscr D}^b(#1)}
\newcommand{\C}[1]{{\mathscr C}(#1)}

\newcommand{\K}[1]{{\mathscr K}(#1)}

\newcommand{\Kb}[1]{{\mathscr K}^b(#1)}
\newcommand{\modcat}{\ensuremath{\mbox{{\rm -mod}}}}
\newcommand{\indmod}{\ensuremath{\mbox{{\rm -ind}}}}

\newcommand{\stmodcat}[1]{#1\mbox{{\rm -{\underline{mod}}}}}

\newcommand{\pmodcat}[1]{#1\mbox{{\rm -proj}}}
\newcommand{\imodcat}[1]{#1\mbox{{\rm -inj}}}

\newcommand{\opp}{^{\rm op}}

\newcommand{\Hom}{{\rm Hom}}

\newcommand{\Ext}{{\rm Ext}}

\newcommand{\lra}{\longrightarrow}
\newcommand{\lraf}[1]{\stackrel{#1}{\lra}}
\newcommand{\ra}{\rightarrow}

\title{ \bf  Extension dimensions: derived equivalences and stable equivalences
\\
\footnotetext{Email addresses: zhangjb@cnu.edu.cn, zhengjunling@cjlu.edu.cn}
}

\author {Jinbi Zhang, Junling Zheng
\thanks{Corresponding author} \\
{\it \scriptsize  School of Mathematical Sciences, Peking University, 100871 Beijing, P. R. China }\\
{\it \scriptsize  Department of Mathematics, China Jiliang University, Hangzhou, 310018, Zhejiang Province, P. R. China}
}
\date{}
\begin{document}



\maketitle
\begin{abstract}
We show that the difference of the extension dimensions of two derived equivalent algebras is bounded above by the minimal length of a tilting complex associated with a derived equivalence, and that the extension dimension is an invariant under the stable equivalence. In addition, we provide two sufficient conditions such that the extension dimension is an invariant under particular derived equivalences.
\end{abstract}

\section{Introduction}
Rouquier introduced the dimension of a triangulated category in \cite{r06,r08}.
It is an invariant that measures how quickly the triangulated category can be built from one object.
This dimension also plays an important role in the representation theory of Artin algebras (see \cite{han2009derived, krause2006rouquier, opp09, r06}). Similar to the dimension of triangulated categories, the extension dimension of abelian categories was introduced by Beligiannis in \cite{b08}. Let $A$ be an Artin algebra. We denote by $\ed(A)$ the extension dimension of the category of all finitely generated left $A$-modules. Beligiannis \cite{b08} proved that $\ed(A)=0$ if and only if $A$ is representation-finite. It means that the extension dimension of an Artin algebra gives a reasonable way of measuring how far an algebra is from being representation-finite. Though many upper bounds have been found for the extension dimension of a given Artin algebra (see \cite{b08,zh22, zheng2020}), the precise value of its extension dimension is very hard to directly compute. One possible strategy is to study the relationship between the extension dimensions of ``nicely" related algebras. Specifically, we study in this paper the following question:

\emph{Suppose two Artin algebras $A$ and $B$ are derived equivalent or stably equivalent, how are the relationships of their extension dimensions?}

As we known, derived equivalences play an important role in the representation theory of Artin algebras and finite groups (see \cite{h88, Xi18}), while the Morita theory of derived categories of rings by Rickard \cite{r1989}  and the Morita theory of derived categories of differential graded algebras by Keller \cite{k94} provide a useful tool to understand homological properties of these equivalent algebras. Many homological invariants of derived equivalences have been discovered, for example Hochschild homology \cite{r91}, cyclic homology \cite{k98}, algebraic $K$-groups \cite{ds04} and the number of non-isomorphic simple modules \cite{r1989}. Though derived equivalences do not always preserve homological dimensions of algebras and modules, they still can provide a useful tool to understand some homological properties of algebras. For example, the differences of global and finitistic dimensions of two derived equivalent algebras are bounded above the length of a tilting complex inducing a derived equivalence (see \cite[Section 12.5(b)]{gr92},\cite{h93},\cite{px09}). We hope to bound the difference of extension dimensions of derived equivalent algebras in terms of lengths of tilting complexes. Recall that the \emph{length} of a radical complex $\cpx{X}$ in $\Kb{A}$ is defined to be
$$\ell(\cpx{X}) = \sup\{s\mid X^s\neq 0\} - \inf\{t\mid X^t\neq 0\} +1.$$
Define the \emph{length} of an arbitrary complex $\cpx{Y}$ in $\Kb{A}$ to be the length of the unique radical complex that is isomorphic to $\cpx{Y}$ in $\Kb{A}$ (Lemma \ref{radical}). One of main results reads as the following theorem.

\begin{Theo}\label{main-der}{\rm (Theorem \ref{der-thm})} Let $F:\Db{A}\lraf{\sim} \Db{B}$ be a derived equivalence between Artin algebras. Then $|\ed(A)-\ed(B)|\le \ell(F(A))-1.$
\end{Theo}

In representation theory, another important equivalence is the stable equivalence of Artin algebras. In \cite{MV1990}, Mart\'inez-Villa proved that stable equivalences preserve the global and dominant dimensions of algebras without nodes.
Recently, Xi-Zhang  \cite{Xi2022} showed that the delooping levels, $\phi$-dimensions and $\psi$-dimensions of Artin algebras are invariants of stable equivalences of algebras with no nodes. Guo \cite{Guo05} showed that stable equivalences preserve the representation dimensions of Artin algebras (this was already proved by Xi in \cite{Xi02} for stable equivalence of Morita type). We will establish the version of extension dimensions under stable equivalences.

\begin{Theo}\label{main-st}{\rm (Theorem \ref{st-thm})}
Let $A$ and $B$ be stably equivalent Artin algebras. Then $\ed(A)=\ed(B)$.
\end{Theo}

In general, one could not hope that the extension dimension is an invariant under derived equivalences. For example, a representation-finite algebra can derived equivalent to a representation-infinite algebra (see Example \ref{ex-silt}). In this article, we provide two sufficient conditions such that derived equivalences preserve extension dimensions.

One well-developed approach is based on the special derived equivalences.
For finite dimensional selfinjective algebras, Rickard \cite{r89} showed that each derived equivalence induces a stable equivalence.
Hu-Xi \cite{hx10} generalized the result of Rickard by introducing a new class of derived equivalences, called almost $\nu$-stable derived equivalences. They proved that every almost $\nu$-stable derived equivalence always induces a stable equivalence. An application of Theorem \ref{main-st} is the following result.

\begin{Koro}\label{main-nu}{\rm (Corollary \ref{nu-s})}
Let $A$ and $B$ be almost $\nu$-stable derived equivalent finite dimensional algebras. Then $\ed(A)=\ed(B)$.
\end{Koro}
The second approach is based on the theory of $2$-term silting complexes, which generalizes classical tilting theory. Hoshino-Kato-Miyachi \cite{HKM02} proved that each $2$-term silting complex $\cpx{P}$ over an algebra $A$ can induce a torsion pair $(\T(\cpx{P}),\F(\cpx{P}))$ in $A\modcat$.
Recently, Buan-Zhou \cite{Buan2016} gave a generalization of the Brenner-Butler tilting theorem (see \cite{bb79,hr82}), called the silting theorem.
This theorem described the relations of torsion pairs between $A\modcat$ and $B\modcat$, where $B=\End_{\Db{A}}(\cpx{P})$.
This provides a basic framework for comparing the extension dimensions of derived equivalent algebras.
Recall that a $2$-term silting complex $P^{\bullet}$ over $A$ is called \emph{separating} if the induced torsion pair $(\T(\cpx{P}),\F(\cpx{P}))$ is split.
Then one of main results can be presented as follows.

\begin{Theo}\label{main-silt}{\rm (Theorem \ref{theo-silt})}
Let $A$ be an Artin algebra, $P^{\bullet}$ a $2$-term separating silting complex over $A$, and $B:=\End_{\Db{A}}(P^{\bullet})$. If $\id(_AX)\le 1$ for each $X\in \F(P^{\bullet})$, then $\ed(A)=\ed(B)$.
\end{Theo}

The paper is outlined as follows: In Section 2, we
recall some basic notations, definitions and facts required in proofs. In Section 3, we compare the extension dimensions of derived equivalent algebras. We first in Section \ref{der-var} gives a proof of Theorem \ref{main-der}. Section \ref{der-inv} then provides a sufficient condition for the extension dimension to be a derived invariant and gives several examples to illustrate the necessity of some assumptions in Theorem \ref{main-silt}. In Section 4, we prove Theorem \ref{main-st} and present an example to illustrate this main result. Corollary \ref{main-nu} is the direct consequence of Theorem \ref{main-st}.

\section{Preliminaries}
In this section, we shall fix some notations, and recall some definitions.

\subsection{Stable equivalences and derived equivalences}
Throughout this paper, $k$ is an arbitrary but fixed commutative Artin ring. Unless stated otherwise, all algebras are Artin $k$-algebras with unit, and all modules are finitely generated unitary left modules; all categories will be $k$-categories
and all functors are $k$-functors.

Let $A$ be an Artin algebra. We denote by $A\modcat$ the category of all finitely generated left $A$-modules, and by $A\indmod$ the set of isomorphism classes of indecomposable finitely generated $A$-modules. All subcategories of $A\modcat$ are full, additive and closed under isomorphisms.
For a class of $A$-modules $\mathcal{X}$, we write $\add(\mathcal{X})$ for the smallest full subcategory of $A\modcat$ containing $\mathcal{X}$ and closed under finite direct sums and direct summands.
When $\mathcal{X}$ consists of only one object $X$, we write $\add(X)$ for $\add(\mathcal{X})$.
Particularly, $\add({_A}A)$ is exactly the category of projective $A$-modules and also denoted by $A\pmodcat$. We denote by $\mathscr{P}_A$ and $\mathscr{I}_A$ the set of isomorphism classes of indecomposable projective and injective $A$-modules, respectively.
Let $X$ be an $A$-module.
If $f:P\ra X$ is the projective cover of $X$ with $P$ projective, then the kernel of $f$ is called the \emph{syzygy} of $X$, denoted by $\Omega(X)$.
Dually, if $g:X\ra I$ is the injective envelope of $X$ with $I$ injective, then the cokernel of $g$ is called the \emph{cosyzygy} of $X$, denoted by $\Omega^{-1}(X)$. Additionally, let $\Omega^0$ be the identity functor in $A\modcat$ and $\Omega^1:=\Omega$. Inductively, for any $n\ge 2$, define $\Omega^n(X):=\Omega^1(\Omega^{n-1}(X))$ and $\Omega^{-n}(X):=\Omega^{-1}(\Omega^{-n+1}(X))$.
We denoted by $\pd(_AX)$ and $\id(_AX)$ the projective and injective dimension, respectively.

Let $A^{\rm {op}}$ be the opposite algebra of $A$, and $\DD:=\Hom_k(-,E(k/\rad(k)))$ the usual duality from $A\modcat$ to $A^{\rm {op}}\modcat$, where $\rad(k)$ denotes the radical of $k$ and $E(k/\rad(k))$ denotes the injective envelope of $k/\rad(k)$. The duality $\Hom_A(-, A)$ from $\pmodcat{A}$ to $\pmodcat{A\opp}$ is denoted by $^*$, namely for each projective $A$-module $P$, the projective $A\opp$-module $\Hom_A(P, A)$ is written as $P^*$. We write $\nu_A$ for the Nakayama functor $\DD\Hom_A(-,A): \pmodcat{A}\ra \imodcat{A}$. An $A$-module $X$ is called a \emph{generator} if $A\in\add(X)$, \emph{cogenerator} if $\DD(A_A)\in\add(X)$, and \emph{generator-cogenerator} if it is both a generator and cogenerator in $A\modcat$.

We denoted by $\stmodcat{A}$ the stable module category of $A$ modulo projective modules. The objects are the same as the objects of $A\modcat$, and the homomorphism set $\underline{\Hom}_A(X,Y)$ between $X$ and $Y$ is given by the quotients of $\Hom_A(X,Y)$  modulo those homomorphisms that factorize through a projective $A$-module. This category is usually called the \emph{stable module category} of $A$. Dually, We denoted by  $A\mbox{-}\overline{\mbox{mod}}$ the stable module category of $A$ modulo injective modules.
Two algebras $A$ and $B$ are said to be \emph{stably equivalent} if the two stable categories $A\stmodcat$ and $B\stmodcat$ are equivalent as additive categories.

\medskip
Let $\cal C$ be an additive category. For two morphisms $f:X\rightarrow Y$ and $g:Y\rightarrow Z$ in $\cal C$, their composition is denoted by $fg$, which is a morphism from $X$ to $Z$. But for two functors $F:\mathcal{C}\ra \mathcal{D}$ and $G:\mathcal{D}\ra\mathcal{E}$ of categories, their composition is written as $GF$.

A complex $\cpx{X}=(X^i, d_X^i)$ over $\mathcal{C}$ is a sequence of objects $X^i$ in
$\cal C$ with morphisms $d_{\cpx{X}}^{i}:X^i\ra X^{i+1}$ such that $d_{\cpx{X}}^{i}d_{\cpx{X}}^{i+1}=0$ for all $i \in {\mathbb Z}$. We denote by $\C{C}$  the category of complexes over $\cal C$ , and by $\K{\mathcal C}$ the homotopy category of complexes over $\mathcal{C}$. If $\cal C$ is an abelian category, then we denote by $\D{\cal C}$ the derived category of complexes over $\cal C$. Let $\Kb{\mathcal C}$
and $\Db{\mathcal C}$ be the full subcategory of $\K{\cal C}$ and $\D{\cal C}$ consisting of
bounded complexes over $\mathcal{C}$, respectively. For a given algebra $A$, we simply write $\C{A}$, $\K{A}$ and $\D{A}$ for $\C{A\modcat}$, $\K{A\modcat}$ and $\D{A\modcat}$, respectively. Similarly, we write $\Kb{A}$ and $\Db{A}$ for $\Kb{A\modcat}$ and $\Db{A\modcat}$, respectively.
It is known that $\K{A}$, $\D{A}$, $\Kb{A}$ and $\Db{A}$ are triangulated categories. For a complex $\cpx{X}$ in $\K{A}$ or $\D{A}$, the complex $\cpx{X}[1]$ is obtained from $\cpx{X}$ by shifting $\cpx{X}$ to the left by one degree.

Let $A$ be an Artin algebra.
A homomorphism $f: X\ra Y$ of $A$-modules is said to be a \emph{radical homomorphism} if, for any module $Z$ and homomorphisms $h: Z\ra X$ and $g: Y\ra Z$, the composition $hfg$ is not an isomorphism.
For a complex $(X^i, d_{\cpx{X}}^i)$ over $A\modcat$, if all $d_{\cpx{X}}^i$ are radical homomorphisms, then it is called a \emph{radical complex}, which has the following properties.
\begin{Lem}\label{radical}
{\rm (\cite[pp. 112-113]{hx10})}
Let $A$ be an Artin algebra.

$(1)$ Every complex over $A\modcat$ is isomorphic to a radical complex in $\K{A}$.

$(2)$ Two radical complexes $\cpx{X}$ and $\cpx{Y}$ are isomorphic in $\K{A}$ if and only if they are isomorphic in $\C{A}$.
\end{Lem}

Two algebras $A$ and $B$ are said to be \emph{derived equivalent} if their derived categories $\Db{A}$ and $\Db{B}$ are equivalent as
triangulated categories. In \cite{r1989}, Rickard proved that $A$ and $B$ are derived equivalent if and only if there exists a bounded complex $\cpx{T}$ of finitely generated projective $A$-modules such that $B\simeq\End_{\Db{A}}(\cpx{T})$ and

(1) $\Hom_{\Db{A}}(\cpx{T},\cpx{T}[i])=0$ for all $i\ne 0$;

(2) $\Kb{\pmodcat{A}}={\rm thick}(T^{\bullet})$, where ${\rm thick}(T^{\bullet})$ is the smallest triangulated subcategory of $\Kb{\pmodcat{A}}$ containing $T^{\bullet}$ and closed under finite direct sums and direct summands.

{\parindent=0pt} A complex in $\Kb{\pmodcat{A}}$ satisfying the above two conditions is called a \emph{tilting complex} over $A$. It is known that, given a derived equivalence $F$ between $A$ and $B$, there is a unique (up to isomorphism) tilting complex $\cpx{T}$ over $A$ such that $F(\cpx{T})\simeq B$. This complex $\cpx{T}$ is called a tilting complex \emph{associated} to $F$.

\begin{Lem}\label{der-lem}{\rm (\cite[Lemma 2.1]{hx10})}
Let $A$ and $B$ be two algebras, and let $F: \Db{A}\lra \Db{B}$ be a derived equivalence with a quasi-inverse $F^{-1}$. Then $F(A)$ is isomorphic in $\Db{B}$ to a complex $\cpx{\bar{T}}\in\Kb{\pmodcat{B}}$  of the form
$$0\lra \bar{T}^0\lra \bar{T}^1\lra\cdots\lra\bar{T}^n\lra 0$$
for some $n\geq 0$ if and only if $F^{-1}(B)$ is isomorphic  in $\Db{A}$ to a complex $\cpx{T}\in\Kb{\pmodcat{A}}$ of the form
$$0\lra T^{-n}\lra\cdots\lra T^{-1}\lra T^0\lra 0.$$
\end{Lem}
A special class of derived equivalences can be constructed by tilting modules. Recall that an $A$-module $T$ is said to be a \emph{tilting module} if $T$ satisfied the following three conditions: (1) $\pd(_AT)\le n$, (2) $\Ext_A^i(T,T)=0$ for all $i>0$, and (3) there exists an exact sequence $0\ra A\ra T_0\ra \cdots\ra T_n\ra 0$ in $A\modcat$ with each $T_i$ in $\add(_AT)$. Let $P^{\bullet}(T): 0\ra  P_{n}\ra P_{n-1}\ra \cdots\ra P_{0}\ra 0$ be a projective resolution of $T$. Clearly, $P^{\bullet}(T)$ is a tilting complex over $A$ and $\End_{A}(T) \simeq \End_{\Db{A}}(P^{\bullet}(T))$ as algebras.

\subsection{Extension dimensions}
In this subsection, we shall recall the definition and some results of the extension dimensions (see \cite{b08,zheng2020}).

Let $A$ be an Artin algebra. We denote by $A\modcat$ the category of all finitely generated left $A$-modules. For a class $\T$ of $A$-modules, we denote by $\add(\T)$ the smallest full subcategory of $A\modcat$ containing $\T$ and closed under finite direct sums and direct summands. When $\T$ consists of only one object $T$, we write $\add(T)$ for $\add(\T)$.
Let $\T_1,\T_2,\cdots,\T_n$ be subcategories of $A\modcat$. Define
\begin{align*}
\T_1\bullet \T_2
:&=\add(\{X\in A\modcat \mid \mbox{there exists an exact sequence } 0\lra T_1\lra  X \lra T_2\lra 0\\
&\qquad\qquad\qquad\qquad\quad\mbox{in } A\modcat \mbox{ with }T_1 \in \T_1\mbox{ and }T_2 \in \T_2\})\\
&=\{X\in A\modcat \mid \mbox{there exists an exact sequence } 0\lra T_1\lra  X\oplus X' \lra T_2\lra 0\\
&\qquad\qquad\qquad\quad\mbox{ in } A\modcat \mbox{ for some } A\mbox{-module } X' \mbox{ with }T_1 \in \T_1\mbox{ and }T_2 \in \T_2\}.
\end{align*}
Inductively, define
$$\T_{1}\bullet  \T_{2}\bullet \dots \bullet\T_{n}:=\add(\{X\in A\modcat \mid \mbox{there exists an exact sequence }
0\lra T_1\lra  X \lra T_2\lra 0$$
$$\mbox{in } A\modcat \mbox{ with }T_1 \in \T_1\mbox{ and }T_2 \in \T_{2}\bullet \dots \bullet\T_{n}\}).$$
Thus $X\in \T_{1}\bullet  \T_{2}\bullet \dots \bullet\T_{n}$ if and only if there exist the following exact sequences
\begin{equation*}
\begin{cases}
\xymatrix@C=1.5em@R=0.1em{
0\ar[r]& T_1 \ar[r]& X\oplus X_1' \ar[r]& X_2 \ar[r]& 0,\\
0\ar[r]& T_2 \ar[r]& X_2\oplus X_2' \ar[r]& X_3 \ar[r]& 0,\\
&&\vdots&&\\
0\ar[r]& T_{n-1} \ar[r]& X_{n-1}\oplus X_{n-1}' \ar[r]& X_n \ar[r]& 0,}
\end{cases}
\end{equation*}
for some $A$-modules $X_i'$ such that $T_i\in \T_i$ and $X_{i+1}\in \T_{i+1}\bullet \T_{i+2}\bullet \cdots \bullet \T_n$ for $1\le i\le n-1$.

For a subcategory $\T$ of $A\modcat$, set $[\T]_{0}:=\{0\}$, $[\T]_{1}:=\add(\T)$,  $[\T]_{n}=[\T]_1\bullet [\T]_{n-1}$ for any $n\ge 2$. If $T\in A\modcat$, we write $[T]_{n}$ instead of $[\{T\}]_{n}$.
\begin{Def}\label{ed-def}{\rm (\cite{b08})} {\rm The extension dimension of $A\modcat$ is defined to be
$$\ed(A):=\inf\{n\ge 0\mid A\modcat=[T]_{n+1}\mbox{ with } T\in A\modcat\}.$$
}
\end{Def}

\begin{Lem} \label{ext-lem}
Let $A$ be an Artin algebra.

{\rm(1) (\cite[Example 1.6]{b08})} $A$ is representation finite if and only if $\ed(A)=0$.

{\rm(2) (\cite[Example 1.6]{b08})} $\ed(A)\le \ell\ell(A)-1$, where $\ell\ell(A)$ stands for the Loewy length of $A$.

{\rm(3) (\cite[Corollary 3.6]{zheng2020})} $\ed(A)\le \gd(A)$.

{\rm(4) (\cite[Corollary 3.15]{zh22})} $\ed(A)\le \ell\ell^{\infty}(A)+\max\{\pd(S)\mid S\mbox{ is simple with }\pd(S)<\infty\}$, where $\ell\ell^{\infty}(A)$ stands for the infinite-layer length of $A$ {\rm(\cite{hlh09, hlh13})}.

\end{Lem}
\begin{Bsp}{\rm
Let $A$ be the Beilinson algebra $kQ/I$ with quiver $Q$
$$\xymatrix{
&0 \ar@/_1pc/[r]_{x_{n}}\ar@/^1pc/[r]^{x_{0}}_{\vdots}
&1\ar@/_1pc/[r]_{x_{n}}\ar@/^1pc/[r]^{x_{0}}_{\vdots}
&2\ar@/_1pc/[r]_{x_{n}}\ar@/^1pc/[r]^{x_{0}}_{\vdots}
&3&\cdots &n-1\ar@/_1pc/[r]_{x_{n}}\ar@/^1pc/[r]^{x_{0}}_{\vdots}&n
}$$
and relations $I=(x_{i}x_{j}-x_{j}x_{i})$ for $0 \le i, j \le n$. By \cite[Example 3.4]{zh22}, we know that $\ed(A)=n$.
We see that the extension dimension may be very large.
}
\end{Bsp}

\begin{Def}\label{wrd-def}{\rm (\cite[Defition 4.5(2)]{Iyama2003})
Let $M$ be an $A$-module. Then the \emph{weak $M$-resolution dimension} of an $A$-module $X$ is defined to be
\begin{align*}
M\mbox{-}\wrd(X):&=\inf \{n\in\mathbb{N} \mid \mbox{ there is an exact sequence } 0 \ra M_{n} \ra M_{n-1}\ra\cdots \ra M_{0}\ra X'\ra 0 \\
&\qquad\qquad\qquad\;\mbox{with all } M_{i}\in \add(M) \mbox{ and }X\in\add(X')\}
\\
&=\inf \{n\in\mathbb{N} \mid \mbox{ there is an exact sequence } 0 \ra M_{n} \ra M_{n-1}\ra\cdots \ra M_{0}\ra X\oplus Y\ra 0 \\
&\qquad\qquad\qquad\mbox{ for some }A\mbox{-module } Y\mbox{ with all } M_{i}\in \add(M)\}.
\end{align*}
Here we set $\inf\varnothing=\infty$; and the \emph{weak $M$-resolution dimension} of algebra $A$ is defined to be
\begin{center}
$M$-$\wrd(A):=\sup\{M$-$\wrd(X) \mid X\in A\modcat\}$;
 \end{center}
the \emph{weak resolution dimension} of algebra $A$ is defined to be
\begin{center}
$\wrd(A):=\inf\{M$-$\wrd(A)\mid M\in A\modcat\}$.
\end{center}

Dually, we can define the \emph{weak $M$-coresolution dimension}, the \emph{weak $M$-coresolution dimension} of $A$ and the \emph{weak coresolution dimension} of $A$. By \cite[Definition 4.5]{Iyama2003}, we know that the weak resolution dimension of $A$ and the weak coresolution dimension of $A$ are the same.
}\end{Def}
Recall that a module is called \emph{basic} if it is a direct sum of non-isomorphic indecomposable modules. Then each module $M$ admits a unique decomposition (up to isomorphism) $M\simeq M_b\oplus M_0$ such that $M_b$ is basic and $M_0\in\add(M_b)$. Then $M\mbox{-}\wrd(A)=M_b\mbox{-}\wrd(A)$ and
\begin{center}
$\wrd(A)=\min\{M$-$\wrd(A)\mid M\in A\modcat \mbox{ and } M \mbox{ is basic} \}$.
\end{center}

In the following lemma, we mention a few basic properties of weak resolution dimensions.
\begin{Lem}\label{weak-prop}
Let $M$ be an $A$-module. Then

$(1)$ For $A$-modules $X_i,1\le i\le n$, we have $$M\mbox{-}\wrd(\bigoplus\limits_{i=1}^nX_i)
=\sup\{M\mbox{-}\wrd(X_i)\mid 1\le i\le n\}.$$

$(2)$ $M\mbox{-}\wrd(A)=\sup\{M\mbox{-}\wrd(X)\mid X\in A\indmod\}$,
where $A\indmod$ stands for the set of isomorphism classes of indecomposable finitely generated $A$-modules.

$(3)$ Fix an $A$-module $M_0$, we have $(M\oplus M_0)\mbox{-}\wrd(A)\le M\mbox{-}\wrd(A)$. In particular, $$\wrd(A)=\inf\{(M'\oplus M_0)\mbox{-}\wrd(A)\mid M'\in A\modcat \mbox{ and } M' \mbox{ is basic}\}.$$
\end{Lem}

\medskip
In \cite{zheng2020}, Zheng-Ma-Huang proved that the weak resolution dimension and the extension dimension of an Artin algebra are coincide. We think the proof of \cite[Theorem 3.5]{zheng2020} is incomplete. Here, we shall give a detailed proof. To prove the result, we need the following lemma.

\begin{Lem}\label{wrd-lem} Let $A$ be an Artin algebra.

{\rm (1) (see \cite[Lemma 4.6]{zheng22})} Given the exact sequence $0 \ra X\ra Y\ra Z\ra 0$
in $A\modcat$, we can get the following exact sequence
$$0 \lra \Omega^{i+1}(Z) \lra \Omega^{i}(X)\oplus P_{i}\lra \Omega^{i}(Y)\lra 0$$
for some projective module $P_{i}$ in $A\modcat$, where $i\ge 0.$

{\rm (2) (\cite[Lemma 3.5]{zh22})} Let $0\ra M_{n}\ra M_{n-1}\ra \cdots\ra M_0\ra X\ra 0$ be an exact sequence in $A\modcat$. Then
$$X\in [M_0]_1\bullet[\Omega^{-1}(M_{1})]_1
\bullet\cdots\bullet
[\Omega^{-n+1}(M_{n-1})]_1\bullet [\Omega^{-n}(M_n)]_1
\subseteq [\bigoplus_{i=0}^n\Omega^{-i}(M_i)]_{n+1}.$$
\end{Lem}

\begin{Lem}\label{wrd=ed}
For an Artin algebra $A$, we have $\wrd(A)=\ed(A)$.
\end{Lem}
{\it Proof.} Suppose $\ed(A)=n$. By Definition \ref{ed-def}, there exists an $A$-module $T$ such that $$A\modcat=[T]_{n+1}.$$
Then, for an $A$-module $X_0$, there are short exact sequences
\begin{equation*}
\begin{cases}
\xymatrix@C=1.5em@R=0.1em{
0\ar[r]& T_0\ar[r] & X_0\oplus X_{0}'\ar[r]& X_{1}\ar[r]&0,\\
0\ar[r]& T_1\ar[r] & X_{1}\oplus X_{1}'\ar[r]& X_{2}\ar[r]&0,\\
&&\vdots&&\\
0\ar[r]& T_{n-1}\ar[r] & X_{n-1}\oplus X_{n-1}'\ar[r]& X_{n}\ar[r]&0.
}
\end{cases}
\end{equation*}
in $A\modcat$ for some $A$-modules $X_i'$ such that $T_i\in [T]_1$ and $X_{i+1}\in [T]_{n-i}$ for each $0\le i\le n-1$. Furthermore, we have short exact sequences
\begin{equation*}
\begin{cases}
\xymatrix@C=1em@R=0.1em{
0\ar[r]& T_0\ar[r]& X_0\oplus (\bigoplus\limits_{j=0}^{n-1} X_j')\ar[r]& X_{1}\oplus (\bigoplus\limits_{j=1}^{n-1} X_j')\ar[r]& 0,\\
0\ar[r]& T_1\ar[r]& X_1\oplus (\bigoplus\limits_{j=1}^{n-1} X_j')\ar[r]& X_{2}\oplus (\bigoplus\limits_{j=2}^{n-1} X_j')\ar[r]& 0,\\
&&\vdots&&\\
0\ar[r]& T_{n-2}\ar[r]& X_{n-2}\oplus (\bigoplus\limits_{j=n-2}^{n-1} X_j')\ar[r]& X_{n-1}\oplus X_{n-1}'\ar[r]& 0,\\
0\ar[r]& T_{n-1}\ar[r]& X_{n-1}\oplus X_{n-1}'\ar[r]& X_{n}\ar[r]& 0.
}
\end{cases}
\end{equation*}
By Lemma \ref{wrd-lem}(1), we have the following short exact sequences
\begin{equation}\label{short-ed}
\begin{cases}
\xymatrix@C=1em@R=0.1em{
0\ar[r]& \Omega^1(X_{1}\oplus (\bigoplus\limits_{j=1}^{n-1} X_j'))\ar[r]& T_0\oplus P_0\ar[r]& X_0\oplus (\bigoplus\limits_{j=0}^{n-1} X_j')\ar[r]& 0,\\
0\ar[r]&\Omega^2(X_{2}\oplus (\bigoplus\limits_{j=2}^{n-1} X_j'))\ar[r]& \Omega^1(T_1\ar[r])\oplus P_1& \Omega^1(X_1\oplus (\bigoplus\limits_{j=1}^{n-1} X_j'))\ar[r]& 0,\\
&&\vdots&&\\
0\ar[r]&\Omega^{n-1}(X_{n-1}\oplus X_{n-1}')\ar[r]& \Omega^{n-2}(T_{n-2})\oplus P_{n-2}\ar[r]& \Omega^{n-2}(X_{n-2}\oplus (\bigoplus\limits_{j=n-2}^{n-1} X_j'))\ar[r]& 0,\\
0\ar[r]& \Omega^n(X_{n})\ar[r]& \Omega^{n-1}(T_{n-1})\oplus P_{n-1}\ar[r]& \Omega^{n-1}(X_{n-1}\oplus X_{n-1}')\ar[r]& 0
}
\end{cases}
\end{equation}
where all $_AP_i$ are projective. By the short exact sequences (\ref{short-ed}), we get the following long short sequence
\begin{align}\label{long-ed}
0\lra \Omega^n(X_n)\lra \Omega^{n-1}(T_{n-1})\oplus P_{n-1}\lra\cdots \lra \Omega^1(T_1)\oplus P_1\lra T_0\oplus P_0\lra X_0\oplus (\bigoplus_{j=i}^{n-1} X_0')\lra 0.
\end{align}
Set $N:=(\bigoplus_{0\le i\le n}\Omega^i(T))\oplus A$. It follows from $X_n\in [T]_1$, $P_i\in \add(_AA)$ and $T_i\in[T]_1$ that $\Omega^n(X_n),P_i,\Omega^{i}(T_{i})\in \add(_AN)$ for $0\le i\le n-1$. By Definition \ref{wrd-def} and the long exact sequence (\ref{long-ed}), we get
$$N\mbox{-}\wrd(X_0)\le n
\quad\mbox{and}\quad N\mbox{-}\wrd(A)\le n.$$
By Definition \ref{wrd-def}, we get $$\wrd(A)\le N\mbox{-}\wrd(A)\le n=\ed(A).$$

Conversely, suppose $\wrd(A)=m$. By Definition \ref{wrd-def}, there exists an $A$-module $M$ such that, for any $A$-module $X$, there is an exact sequence
$$0\lra M_{m}\lra M_{m-1}\lra \cdots\lra M_0\lra X\oplus Y\lra 0$$
in $A\modcat$ for some $A$-module $Y$ such that $M_i\in \add(M)$ for $0\le i\le m$. By Lemma \ref{wrd-lem}(2), we have $X\oplus Y\in [\bigoplus_{i=0}^m\Omega^{-i}(M_i)]_{m+1}$ and $X\in [\bigoplus_{i=0}^m\Omega^{-i}(M_i)]_{m+1}$. Then $A\modcat= [\bigoplus_{i=0}^m\Omega^{-i}(M_i)]_{m+1}$. By Definition \ref{ed-def}, $$\ed(A)\le m=\wrd(A).$$

Thus we have $\wrd(A)=\ed(A)$.
$\square$

\section{Derived equivalences}
In this section, we discuss the relationships of the extension dimensions of two derived equivalent algebras. In the first subsection, we get how much extension dimensions can vary under derived equivalences. The second subsection provides a sufficient condition such that two derived equivalent algebras have the same extension dimensions.

\subsection{Variance of extension dimensions under derived equivalences}\label{der-var}
In this subsection, we first review some of the basic facts and conclusions, as detailed in reference \cite{hp17}.

\begin{Def} A derived equivalence $F:\Db{A}\ra \Db{B}$ is called \emph{nonnegative} if

$(1)$ $F(X)$ is isomorphic to a complex with zero homology in all negative degrees for all $X\in A\modcat$; and
$(2)$ $F(P)$ is isomorphic to a complex in $\Kb{\pmodcat{A}}$ with zero terms in all negative degrees for all $P\in\pmodcat{A}$.
\end{Def}
\begin{Lem}\label{lem-nonneg}
{\rm (\cite[Lemma 4.2]{hp17})} A derived equivalence $F$ is nonnegative if and only if the tilting complex associated to $F$ is isomorphic in $\Kb{\pmodcat{B}}$ to a complex with zero terms in all positive degrees.
In particular, $F[i]$ is nonnegative for sufficiently small $i$.
\end{Lem}
For every nonnegative derived equivalence $F$, Hu-Xi (see \cite[Section 3]{hx10} or \cite[Section 4]{hp17}) construct a functor  $\overline{F}:\stmodcat{A}\ra \stmodcat{B}$, which is called the \emph{stable functor} of $F$.
This stable functor has the following properties.

\begin{Lem}\label{lem-non}{\rm (\cite[Section 4]{hp17})}
$(1)$ Let $i$ be a nonnegative integer. Then $i$-th syzygy functor $\Omega^{i}_{A}:\stmodcat{A}\ra \stmodcat{A}$ is a stable functor of the derived equivalence $[-i]:\Db{A}\ra \Db{A}$, that is, $\overline{[-i]}\simeq \Omega^i_A$ as additive functors. Particularly, the stable functor of identity functor on $\Db{A}$ is isomorphic to the identity functor on $\stmodcat{A}$.

$(2)$ Let $F:\Db{A}\ra \Db{B}$ and $G:\Db{B}\ra \Db{C}$ be two nonnegative derived equivalences. Then the functors $\overline{G}\circ \overline{F}$ and $\overline{GF}$ are isomorphic.

$(3)$ Let $F:\Db{A}\ra \Db{B}$ be a nonnegative derived equivalence. Suppose that $$0 \lra X\lra Y \lra Z\lra 0$$ is an exact sequence in $A\modcat$. Then there is an exact sequence
$$0 \lra \overline{F}(X) \lra \overline{F}(Y)\oplus Q \lra \overline{F}(Z)\oplus Q'\lra 0$$ in $B\modcat$ for some projective modules $Q$ and $Q'$.
\end{Lem}

\begin{Theo}\label{der-thm} Let $F:\Db{A}\lraf{\sim} \Db{B}$ be a derived equivalence between Artin algebras. Then $|\ed(A)-\ed(B)|\le \ell(F(A))-1.$
\end{Theo}
{\it Proof.} Set $n:=\ell(F(A))-1$. Let $\cpx{T}$ and $\cpx{\bar{T}}$ be the radical tilting complexes associated to $F$ and the quasi-inverse $G$ of $F$, respectively. By applying the shift functor, we can assume that $\cpx{T}\in\Kb{\pmodcat{A}}$ is of the form
$$0\lra T^{-n}\lra\cdots\lra T^{-1}\lra T^0\lra 0.$$
By Lemma \ref{der-lem}, $\cpx{\bar{T}}\in\Kb{\pmodcat{B}}$ is of the form
$$0\lra \bar{T}^0\lra \bar{T}^1\lra\cdots\lra\bar{T}^n\lra 0.$$
By Lemma \ref{lem-nonneg}, we know that $F$ and $G[-n]$ are nonnegative.

By Lemma \ref{wrd=ed}, we know that the extension dimension and the weak resolution dimension of an algebra are the same. We have to show $|\wrd(A)-\wrd(B)|\le n$. We first show $\wrd(B)\le \wrd(A)+ n$. Set $t:=M$-$\wrd(A)$. By Definition \ref{wrd-def}, there exists an $A$-module $M$ such that $\wrd(A)=M$-$\wrd(A)=t$.
Next, we prove $$(\overline{F}(M)\oplus B\oplus D(B_B))\mbox{-}\wrd(B)\le t+n.$$ Indeed, for $Y\in B\modcat$, $\overline{G[-n]}(Y)\in A\modcat$.
By Definition \ref{wrd-def}, we have the following exact sequence
$$0 \lra M_{t}\lraf{f_t} M_{t-1}\lraf{f_{t-1}} \cdots \lraf{f_2} M_{1} \lraf{f_1} M_{0}\lraf{f_0} \overline{G[-n]}(Y)\oplus Z\lra 0$$
in $A\modcat$ for some $A$-module $Z$ with $M_i\in \add(_AM)$ for each $i$. Then there are short exact sequences
\begin{equation}\label{kmk-der}
\begin{cases}
\xymatrix@C=1em@R=0.1em{
0\ar[r]& K_{1} \ar[r]& M_0\ar[r]& \overline{G[-n]}(Y)\oplus Z \ar[r]& 0,\\
0\ar[r]& K_{2} \ar[r]& M_1\ar[r]& K_1 \ar[r]& 0,\\
&&\vdots&&\\
0\ar[r]& K_{t-1} \ar[r]& M_{t-2}\ar[r]& K_{t-2} \ar[r]& 0,\\
0\ar[r]& K_{t} \ar[r]& M_{t-1}\ar[r]& K_{t-1} \ar[r]& 0,}
\end{cases}
\end{equation}
where $K_{j+1}$ is the kernel of $f_j$ for $0\le j\le t-1$. By Lemma \ref{lem-non}(3), applying the stable functor $\overline{F}$ of $F$ to the short exact sequences (\ref{kmk-der}) in $A\modcat$, we can get the following short exact sequences
\begin{equation}\label{fkmk-der}
\begin{cases}
\xymatrix@C=1em@R=0.1em{
0\ar[r]& \overline{F}(K_{1}) \ar[r]& \overline{F}(M_0)\oplus Q_0\ar[r]& \overline{F}(\overline{G[-n]}(Y)\oplus Z )\oplus Q_0' \ar[r]& 0,\\
0\ar[r]& \overline{F}(K_{2}) \ar[r]& \overline{F}(M_1)\oplus Q_1\ar[r]& \overline{F}(K_1)\oplus Q_1' \ar[r]& 0,\\
&&\vdots&&\\
0\ar[r]& \overline{F}(K_{t-1}) \ar[r]& \overline{F}(M_{t-2})\oplus Q_{t-2}\ar[r]& \overline{F}(K_{t-2})\oplus Q_{t-2}'\ar[r]& 0,\\
0\ar[r]& \overline{F}(K_{t}) \ar[r]& \overline{F}(M_{t-1})\oplus Q_{t-1}\ar[r]& \overline{F}(K_{t-1})\oplus Q_{t-1}'\ar[r]& 0,}
\end{cases}
\end{equation}
in $B\modcat$ for some projective $B$-modules $Q_i$ and $Q'_i$ for $0\le i\le t-1$. Combine the short exact sequences (\ref{fkmk-der}), we get the following long exact sequence
\begin{align}\label{long-f-der}
0 \lra \overline{F}(M_{t})\lra \overline{F}(M_{t-1})\oplus Q_{t-1}\lra \cdots\lra
\overline{F}(M_{0})\oplus Q_{0}\oplus Q_{1}'\lra \overline{F}(\overline{G[-n]}(Y)\oplus Z)\oplus Q_{0}'\lra 0
\end{align}
in $B\modcat$. On the other hand, we have the following isomorphisms in $\stmodcat{B}$
\begin{align*}
\overline{F}(\overline{G[-n]}(Y)\oplus Z)
&\simeq \overline{F}(\overline{G[-n]}(Y))\oplus \overline{F}(Z)\\
&\simeq (\overline{F \circ G}\circ \overline{[-n]}(Y))\oplus \overline{F}(Z)
\quad \mbox{(by Lemma \ref{lem-non}(2))}\\
&\simeq  \overline{[-n]}(Y)\oplus \overline{F}(Z)
\quad \mbox{(by Lemma \ref{lem-non}(1))}\\
&\simeq  \Omega^{n}_{B}(Y)\oplus \overline{F}(Z)
\quad \mbox{(by Lemma \ref{lem-non}(1))}.
\end{align*}
By \cite[Theorem 2.2]{Heller60}, there are projective $B$-modules $Q$ and $Q'$ such that $$\overline{F}(\overline{G[-n]}(Y)\oplus Z)\oplus Q'\simeq \Omega^{n}_{B}(Y)\oplus \overline{F}(Z)\oplus Q$$ as $B$-modules. Note that there exist the following two exact sequences
\begin{align}\label{long-y-der} 0 \lra\Omega^{n}_{B}(Y) \lra P_{n-1 }\lra P_{n-2}
\lra\cdots \lra P_{1} \lra P_{0} \lra Y \lra 0,\;\mbox{and}
\end{align}
\begin{align}\label{long-z-der}
0 \lra \overline{F}(Z) \lra J^0\lra J^1
\lra\cdots \lra J^{n-2} \lra J^{n-1} \lra \Omega^{-1}(\overline{F}(Z))\lra 0
\end{align}
in $B\modcat$, where $P_{j}\in \add(_BB)$ and all $J^{j}\in \add(_BD(B_B))$ for $0\le j\le n-1$. By the long exact sequences (\ref{long-f-der}), (\ref{long-y-der}), and (\ref{long-z-der}), we get the following long exact sequence
\begin{align}\label{long-der}\xymatrix{
0\ar[r]& \overline{F}(M_{t})\ar[r]
& \overline{F}(M_{t-1})\oplus Q_{t-1}\ra
 \cdots\ra
 \overline{F}(M_{1})\oplus Q_{1}\oplus Q_{2}' \ar[r]
& \overline{F}(M_{0})\oplus Q_{0}\oplus Q_{1}'\oplus Q'
\ar `r[d] `[l] `[llld] `[dll] [dll]\\
& P_{ n-1 }\oplus J^0\oplus Q \ar[r]
& P_{n-2}\oplus J^1 \ra
\cdots\ra
P_{0}\oplus J^{n-1}\ar[r]
& Y\oplus \Omega^{-1}(\overline{F}(Z))\lra 0.
}\end{align}
It follows from $M_i\in \add(_AM)$ that $\overline{F}(M_i)\in \add(_B\overline{F}(M))$ for $0\le i\le t$. Also $Q_{i},Q_{i}'\in\add(_BB)$ for $0 \le i \le t-1$, $P_j\in \add(_BB)$ and $J^j\in \add(\DD(B_B))$ for $0\le j\le n-1$. By Definition \ref{wrd-def} and the long exact sequence (\ref{long-der}), we get
$$(\overline{F}(M)\oplus B\oplus \DD(B_B))\mbox{-}\wrd(Y)\le t+n
\quad\mbox{and}\quad
\wrd(B)\le t+n=\wrd(A)+n.$$ Similarly, we have $\wrd(A)\le \wrd(B)+n$. Thus $|\wrd(A)-\wrd(B)|\le n.$ By Lemma \ref{wrd=ed}, we get $|\ed(A)-\ed(B)|\le n$.
$\square$

As an immediate consequence of Theorem \ref{der-thm}, we have
\begin{Koro}
Let $A$ be an Artin algebra, $T$ be a tilting $A$-module and $B=\End_{A}(T)$, Then we have
$|\ed(A)-\ed(B)|\le \pd(_AT).$
\end{Koro}
{\it Proof.} Let $n:=\pd(_AT).$ We have the following  minimal projective resolution of $_AT$
$$0\lra  P_{n}\lra P_{n-1}\lra \cdots\lra P_{0}\lra T\lra 0.$$ Moreover, the following complex
$$P^{\bullet}(T):  0\lra  P_{n}\lra P_{n-1}\lra \cdots\lra P_{0}\lra 0 $$
is a tilting complex in $\Db{A}$ and $B\simeq \End_{\Db{A}}(P^{\bullet}(T))$. Thus, by Theorem \ref{der-thm}, we get the result.
$\square$

\subsection{The extension invariants of $2$-term silting complexes}\label{der-inv}
Let $A$ be an Artin algebra. Recall that a pair $(\U,\V)$ of full subcategories of $A\modcat$ is called a \emph{torsion pair} \cite{d66} if the following conditions are satisfied:

(1) $\Hom_A(\U,V)=0$ if and only if $V\in\V$;

(2) $\Hom_A(U,\V)=0$ if and only if $U\in \U$.

The subcategory $\U$ is called the \emph{torsion class} and the subcategory $\V$ is called the \emph{torsion-free class}. It is
known (see \cite[Proposition 1.1]{Assem1990}) that a subcategory $\U$ (respectively, $\V$) of $A$-modules is a torsion class (respectively, torsion-free class) of a torsion pair in $A\modcat$ if and only if $\U$ (respectively, $\V$) is closed under images (respectively, submodules), direct sums and extensions. A torsion pair $(\U,\V)$ is called {\em split} (or sometimes splitting) if each indecomposable $A$-module lies either in $\U$ or in $\V$.

Recall that a complex $P^{\bullet}=(P^i)$ is $2$-term if $P^i=0$ for $i\neq -1,0$. A $2$-term complex $P^{\bullet}\in \Kb{\pmodcat{A}}$ is said to be \emph{silting} \cite{kv88} if $\Hom_{\Kb{\pmodcat{A}}}(P^{\bullet},P^{\bullet}[1])=0$ and ${\rm thick}(P^{\bullet})=\Kb{\pmodcat{A}}$, where ${\rm thick}(P^{\bullet})$ is the smallest triangulated subcategory of $\Kb{\pmodcat{A}}$ containing $P^{\bullet}$ and closed under finite direct sums and direct summands. If, in addition, $\Hom_{\Kb{\pmodcat{A}}}(P^{\bullet},P^{\bullet}[-1])=0$, then it is easy to see that $P^{\bullet}$ is {\em tilting}.

Let $P^{\bullet}$ be a 2-term silting complex in $\Kb{\pmodcat{A}}$, and consider the following two full subcategories of $A\modcat$ given by
\begin{align*}
\T(P^{\bullet}) & =\{U\in A\modcat \mid \Hom_{\Db{A}}(P^{\bullet},U[1])=0\} ,\mbox{ and}\\
\F(P^{\bullet}) & =\{V\in A\modcat \mid \Hom_{\Db{A}}(P^{\bullet},V)=0\}.
\end{align*}

The following lemma is a generalization of the Brenner-Butler tilting theorem (see \cite{bb79,hr82}) to $2$-term silting complexes.
\begin{Lem}\label{equi}{\rm (\cite{Buan2016})}
Let $P^{\bullet}$ be a 2-term silting complex in $\Kb{\pmodcat{A}}$, and let $B=\End_{\Db{A}}(P^{\bullet})$.

$(1)$ Let $\CC(P^{\bullet}):=\{W^{\bullet}\in \Db{A}\mid \Hom_{\Db{A}}(P^{\bullet},W^{\bullet}[i])=0, \;\forall i\neq 0\}$. Then $\CC(P^{\bullet})$ is an abelian category and the short exact sequences in $\CC(P^{\bullet})$ are precisely the triangles in $\Db{A}$ all of whose vertices are objects in $\CC(P^{\bullet})$.

$(2)$ The pairs $(\T(P^{\bullet}),\F(P^{\bullet})$ and $(\F(P^{\bullet})[1],\T(P^{\bullet}))$ are torsion pairs in $A\modcat$ and $\CC(P^{\bullet})$, respectively.

$(3)$ $\Hom_{\Db{A}}(P^{\bullet},-):\CC(P^{\bullet})\ra B\modcat$ is an equivalence of abelian categories. In particular, if $P^{\bullet}$ is a tilting complex, then there exists an equivalence of triangulated categories $F:\Db{A}\ra \Db{B}$ such that $F^{-1}(B)\simeq P^{\bullet}$ and $F(W^{\bullet})\simeq \Hom_{\Db{A}}(P^{\bullet}, W^{\bullet})$ for any $W^{\bullet}\in \CC(P^{\bullet})$.

$(4)$ There is a triangle
$$A\lra P^{\bullet}_1\lraf{f} P^{\bullet}_0\lra A[1]$$
in $\Kb{\pmodcat{A}}$ with $P^{\bullet}_1$, $P^{\bullet}_0\in \add(P^{\bullet})$.

Consider the following $2$-term complex $Q^{\bullet}$ in $\Kb{\pmodcat{B}}$.
$$\xymatrix{0\ar[r]& \Hom_{\Db{A}}(P^{\bullet}, P^{\bullet}_1)\ar[rr]^{\Hom_{\Db{A}}(P^{\bullet}, f)} &&\Hom_{\Db{A}}(P^{\bullet}, P^{\bullet}_0)\ar[r] &0.}$$

$(5)$ $Q^{\bullet}$ is a $2$-term silting complex in $\Kb{\pmodcat{B}}$.

$(6)$ There is an algebra epimorphism $\Phi:A\ra \bar{A}:=\End_{\Db{B}}(Q^{\bullet})$. Moreover, $\Phi$ is an isomorphism if and only if $P^{\bullet}$ is tilting.

$(7)$ Let $\Phi_{\ast}:\bar{A}\modcat\ra A\modcat$ be the inclusion functor induced by $\Phi$ in $(6)$. Then the functor $\Hom_{\Db{A}}(P^{\bullet}, -):\T(P^{\bullet})\ra \F(Q^{\bullet})$ is an equivalence with quasi-inverse $\Phi_{\ast} \Hom_{\Db{B}}(Q^{\bullet}, -[1])$; the functor $\Hom_{\Db{A}}(P^{\bullet}, -[1]):\F(P^{\bullet})\ra \T(Q^{\bullet})$ is an equivalence with quasi-inverse $\Phi_{\ast} \Hom_{\Db{B}}(Q^{\bullet}, -)$.
$$\xymatrix@C=3.2cm{
\T(P^{\bullet})\ar@<1ex>[r]^{\Hom_{\Db{A}}(P^{\bullet}, -)}
&\F(Q^{\bullet})\ar@<1ex>[l]^{\Phi_{\ast} \Hom_{\Db{B}}(Q^{\bullet}, -[1])}
}
\quad\mbox{and}\quad
\xymatrix@C=3.2cm{
\F(P^{\bullet})\ar@<1ex>[r]^{\Hom_{\Db{A}}(P^{\bullet}, -[1])}
&\T(Q^{\bullet}).\ar@<1ex>[l]^{\Phi_{\ast} \Hom_{\Db{B}}(Q^{\bullet}, -)}
}$$

\end{Lem}
{\it Proof.} Note that (1)-(3) is from \cite{HKM02}, (4) is from \cite[Theorem 3.5]{wei13}, (5) and (6) could have been deduced from \cite[Propositions A.3 and A.5]{by13}, and (7) is from \cite[Theorem 1.1]{Buan2016}. $\square$

In the following, the symbol $Q^{\bullet}$ always denotes the induced complex $Q^{\bullet}$. It is a  $2$-term silting complex in $\Kb{\pmodcat{B}}$.

\begin{Def}\label{sep-split}{\rm (\cite[Definition 5.4]{Buan2016})
Let $P^{\bullet}$ be a $2$-term silting complex in $\Kb{\pmodcat{A}}$.

$(1)$ $P^{\bullet}$ is called \emph{separating} if the induced torsion pair ($\T(P^{\bullet})$, $\F(P^{\bullet})$) in $A\modcat$ is split.

$(2)$ $P^{\bullet}$ is called \emph{splitting} if the induced torsion pair ($\T(Q^{\bullet})$, $\F(Q^{\bullet})$) in $B\modcat$ is split.
}\end{Def}

\begin{Lem}\label{sep}
Let $P^{\bullet}$ be a $2$-term silting complex in $\Kb{\pmodcat{A}}$.

$(1)$ {\rm (\cite[Lemma 5.5]{Buan2016})} $P^{\bullet}$ is splitting if and only if $\Ext^{2}_{A}(\T(P^{\bullet}),\F(P^{\bullet}))=0$.

$(2)$ {\rm (\cite[Proposition 5.7]{Buan2016})} If $P^{\bullet}$ is separating, then $P^{\bullet}$ is a tilting complex.

$(3)$ {\rm (\cite[Proposition 3.5]{Hu2021})} Suppose $\id(_AX)\le 1$ for each $X\in \F(P^{\bullet})$. Then $P^{\bullet}$ is separating if and only if $\pd(_BY)\le 1$ for each $Y\in \T(Q^{\bullet})$.
\end{Lem}
\begin{Theo}\label{theo-silt}
Let $A$ be an Artin algebra, $P^{\bullet}$ a $2$-term separating silting complex, and $B:=\End_{\Db{A}}(P^{\bullet})$. If $\id(_AX)\le 1$ for each $X\in \F(P^{\bullet})$, then $\ed(A)=\ed(B)$.
\end{Theo}
{\it Proof.} It follows from Lemma \ref{wrd=ed} that the extension dimension and the weak resolution dimension of an algebra are the same. Thus we have to show that $\wrd(A)=\wrd(B)$.

It follows from $\id(_AX)\le 1$ for each $X\in \F(P^{\bullet})$ that $\Ext^{2}_{A}(\T(P^{\bullet}),\F(P^{\bullet}))=0$. By Lemma \ref{sep}(1), we get $P^{\bullet}$ is splitting, that is, the torsion pair ($\T(Q^{\bullet})$, $\F(Q^{\bullet})$) in $B\modcat$ is split. By Definition \ref{sep-split}(1), since $P^{\bullet}$ is separating, we know that the torsion pair $(\T(P^{\bullet}),\F(P^{\bullet}))$ in $A\modcat$ is split. Denote
$$\H:=\Hom_{\Db{A}}(P^{\bullet},-) \quad\mbox{and}\quad
\E:=\Hom_{\Db{A}}(P^{\bullet},-[1]).$$
By Lemma \ref{equi}(7), we have
$$\H: \T(P^{\bullet})\lraf{\simeq} \F(Q^{\bullet})
\mbox{ and }
\E:\F(P^{\bullet})\lraf{\simeq} \T(Q^{\bullet})$$ as additive categories.

We first prove that $\wrd(A)=0$ if and only if $\wrd(B)=0$. Suppose $\wrd(A)=0$.
By Lemmas \ref{ext-lem}(1) and \ref{wrd=ed}, the weak resolution dimension of an algebra is equal to $0$ if and only if it is representation-finite. Thus $A$ is representation-finite, that is, the number of non-isomorphic indecomposable $A$-modules is finite.
It follows from $\H: \T(P^{\bullet})\lraf{\simeq} \F(Q^{\bullet})$ and $\E:\F(P^{\bullet})\lraf{\simeq} \T(Q^{\bullet})$ as additive categories that the number of non-isomorphic indecomposable $B$-modules in $\F(Q^{\bullet})$ is equal to the number of non-isomorphic indecomposable $A$-modules in $\T(P^{\bullet})$; the number of non-isomorphic indecomposable $B$-modules in $\T(Q^{\bullet})$ is equal to the number of non-isomorphic indecomposable $A$-modules in $\F(P^{\bullet})$.
Thus the number of non-isomorphic indecomposable $B$-modules in $\F(Q^{\bullet})$ or $\T(Q^{\bullet})$ is finite.
Note that the torsion pair $(\T(Q^{\bullet}),\F(Q^{\bullet}))$ in $B\modcat$ is split, namely each indecomposable $B$-module lies either in $\T(Q^{\bullet})$ or in $\F(Q^{\bullet})$.
Thus the number of non-isomorphic indecomposable $B$-modules is finite, that is, $B$ is representation-finite. By Lemmas \ref{ext-lem}(1) and \ref{wrd=ed}, $\wrd(B)=0$. Similarly, if $\wrd(B)=0$, then $\wrd(A)=0$. Thus $\wrd(A)=0$ if and only if $\wrd(B)=0$.

Next, suppose $\wrd(A)\ge 1$ and $\wrd(B)\ge 1$. Set $m:=\wrd(A)$. By Definition \ref{wrd-def}, there exists an $A$-module $M$ such that $\wrd(A)=M$-$\wrd(A)$. Define $N:=B\oplus \H(M)$. We shall show $N$-$\wrd(B)\le m$.
By Lemma \ref{weak-prop}, $$N\mbox{-}\wrd(B)=\sup\{N\mbox{-}\wrd(_BY)\mid Y\in B\indmod\},$$
where $B\indmod$ stands for the set of isomorphism classes of indecomposable finitely generated $B$-modules.
Note that the torsion pair ($\T(Q^{\bullet}),\F(Q^{\bullet})$) in $B\modcat$ is split, that is, each indecomposable $B$-modules lies either in $\T(Q^{\bullet})$ or in $\F(Q^{\bullet})$.
Thus we have to show that $N\mbox{-}\wrd(_BY)\le m$ for each $Y\in \T(Q^{\bullet})$ or $\F(Q^{\bullet})$.

Indeed, if $Y_0\in \F(Q^{\bullet})$, then it follows from $\H: \T(P^{\bullet})\lraf{\simeq} \F(Q^{\bullet})$ that there exists $U\in \T(P^{\bullet})$ such that $\H(U)\simeq Y_0$.
Also $M$-$\wrd(A)=m$.
By Definition \ref{wrd-def}, there is an exact sequence
$$0\lra M_m\lra M_{m-1}\lra \cdots\lra M_0\lra U\oplus U'\lra 0$$
in $A\modcat$ for some $A$-module $U'$ such that all $M_i\in\add(_AM)$.
Equivalently, there is an exact sequence
\begin{align}\label{kmk-silt}
0\ra K_{i+1}\ra M_i\lraf{f_i} K_i\ra 0
\end{align} in $A\modcat$ such that $M_i\in\add(_AM)$ for each $0\le i\le m-1$, where $K_0:=U\oplus U'$ and $K_n:=M_n$.
Since $\cpx{P}$ is $2$-term complex in $\Kb{\pmodcat{A}}$, we have, for any $W\in A\modcat$, $\Hom_{\Db{A}}(P^{\bullet},W[j])\simeq \Hom_{\Kb{\pmodcat{A}}}(P^{\bullet},W[j])=0$ for any $j\neq -1,0$.
For $0\le i\le m-1$, applying the functor $\H$ to the sequence (\ref{kmk-silt}), we obtain the following exact sequence in $B\modcat$:
\begin{align}\label{kmkemk}
0\lra \H(K_{i+1})\lra\H(M_i)\lraf{\H(f_i)}\H(K_i)\lra\E(E_{i+1})\lra\E(M_i)\lra\E(K_i)\lra 0.
\end{align}
Let $C_i$ be the cokernel of $\H(f_i)$. By the long exact sequence (\ref{kmkemk}), we have an exact sequence $$0\lra C_i\lra \E(E_{i+1})\lra\E(M_i)\lra\E(K_i)\lra 0.$$
By Lemma \ref{sep}(3), $\pd(_BZ)\le 1$ if $Z\in \T(Q^{\bullet})$. Note that $\E(E_{i+1}),\E(M_i),\E(K_i)$ are in $\T(Q^{\bullet})$, and have projective dimension no more than $1$. Thus $\pd(_{B}C_i)\le 1$. Let $0\ra Q_{i,1}\ra Q_{i,0}\ra C_i\ra 0$ be a projective resolution of $C_i$ with $Q_{i,j}\in \pmodcat{B}$. Similar to the proof of the horseshoe lemma, we have the following commutative diagram
$$\xymatrix{
 &0\ar[d]&0\ar[d]&0\ar[d]&\\
0\ar[r] &\H(K_{i+1})\ar[r]\ar[d]&\H(K_{i+1})\oplus Q_{i,1}\ar[r]\ar[d]&Q_{i,1}\ar[r]\ar[d]&0\\
0\ar[r] &\H(M_i)\ar[r]\ar[d]&\H(M_i)\oplus Q_{i,0}\ar[r]\ar[d]&Q_{i,0}\ar[r]\ar[d]&0\\
0\ar[r] &{\rm Im}\H(f_i)\ar[r]\ar[d]&\H(K_i)\ar[r]\ar[d]&C_i\ar[r]\ar[d]&0,\\
 &0&0&0&
}$$
where ${\rm Im}\H(f_i)$ is the image of $\H(f_i)$. Then we have a short exact sequence
\begin{align}\label{hqhqk}
0\ra \H(K_{i+1})\oplus Q_{i,1}\ra \H(M_i)\oplus Q_{i,0}\ra \H(K_i)\ra 0.
\end{align}
Combine these short exact sequences (\ref{hqhqk}) for $0\le i\le m-1$, we get a long exact sequence
\begin{align}\label{long-silt}
0\lra \H(M_m)\oplus Q_{m-1,1}\lra\cdots\lra\H(M_1)\oplus Q_{1,0}\oplus Q_{0,1}\lra \H(M_0)\oplus Q_{0,0}\lra \H(U)\oplus\H(U')\lra 0.
\end{align}
Note $N=B\oplus \H(M)$. It follows from $M_i\in \add(_AM)$ and $Q_{ij}\in \add(_BB)$ that $\H(M_i),Q_{ij}\in \add(_BN)$ for $0\le i\le m-1, j=0,1$. By Definition \ref{wrd-def} and the long exact sequence (\ref{long-silt}), we have $$N\mbox{-}\wrd(_B\H(U))\le m.$$
It follows from $Y_0\simeq \H(U)$ as $B$-modules that $$N\mbox{-}\wrd(_BY_0)= N\mbox{-}\wrd(_B\H(U))\le m.$$
Thus we obtain $$N\mbox{-}\wrd(_BY)\le m,\;\forall\;Y\in \F(Q^{\bullet}).$$
On the other hand, by Lemma \ref{sep}(3), we get $\pd(_BY)\le 1$ for each $Y\in \T(Q^{\bullet})$. Also $B\in\add(_BN)$.
Then $$N\mbox{-}\wrd(_BY)\le \pd(_BY)\le 1\le m,\;\forall \;Y\in \T(Q^{\bullet}).$$
Thus we obtain $N\mbox{-}\wrd(_BY)\le m$ for each $Y\in \T(Q^{\bullet})$ or $\F(Q^{\bullet})$. Then
$$N\mbox{-}\wrd(B)\le m.$$ By Definition \ref{wrd-def}, we have $\wrd(B)\le N$-$\wrd(B)$. Then $$\wrd(B)\le m=\wrd(A).$$ Dually, we can prove $\wrd(A)\le \wrd(B)$. Thus $\wrd(A)=\wrd(B)$.
Finally, we have $\ed(A)=\ed(B)$ by Lemma \ref{wrd=ed}.
$\square$

The following corollary is an immediate consequence of the above theorem.
\begin{Koro} Let $A$ be an Artin algebra.

$(1)$ Suppose that $T$ is a separating and splitting tilting $A$-module with $\pd(_AT)\le 1$. Then $\ed(A)=\ed(\End_A(T))$.

$(2)$ Suppose that $P^{\bullet}$ is a $2$-term separating and splitting silting complex. Then $\ed(A/I) =\ed(\End_A(H^0(P^{\bullet})))$, where $I:={\rm ann}_A(H^0(P^{\bullet}))$ is the annihilator of $H^0(P^{\bullet})$.
\end{Koro}
{\it Proof.} $(1)$ Let $P^{\bullet}$ be the minimal projective resolution of $T$. Clearly, $\End_{\Db{A}}(P^{\bullet})\simeq \End_A(T)$ as algebras and $P^{\bullet}$ is a $2$-term silting complex with $$\T(P^{\bullet})=\{U\in A\modcat\mid \Ext_A^1(T,U)=0\}\mbox{ and }\F(P^{\bullet})=\{V\in A\modcat\mid \Hom_A(T,V)=0\}.$$ By \cite[Theorem 3.6, p. 49]{Assem1990}, we have $\id(_AV)=1$ for all $V\in \F(P^{\bullet})$. Then $$\ed(A)=\ed(\End_{\Db{A}}(P^{\bullet}))=\ed(\End_A(T))$$ follows from Theorem \ref{theo-silt}.

(2) By \cite[Proposition 5.1]{Hu2021}, $H^0(P^{\bullet})$ is a separating and splitting tilting $A/I$-module. Then the statement in (2) follows from (1).
$\square$

The following example illustrates that the assumptions in Theorem \ref{theo-silt} are necessary.
\begin{Bsp}\label{ex-silt}{\rm
Let $A$ be an algebra over a field $k$ given by the following quiver $Q_A$
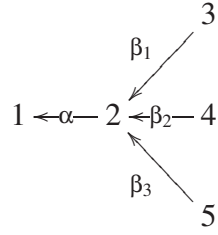
$$\xymatrix{
&&3\ar[dl]_{\beta_1}\\
1&2\ar[l]|-{\alpha}&4\ar[l]|-{\beta_2}\\
&&5\ar[ul]^{\beta_3}
}$$
with relations $\alpha\beta_i=0,\;1\le i\le 3$. The Aulander-Reiten quiver of $A\modcat$ is as follows:
$$
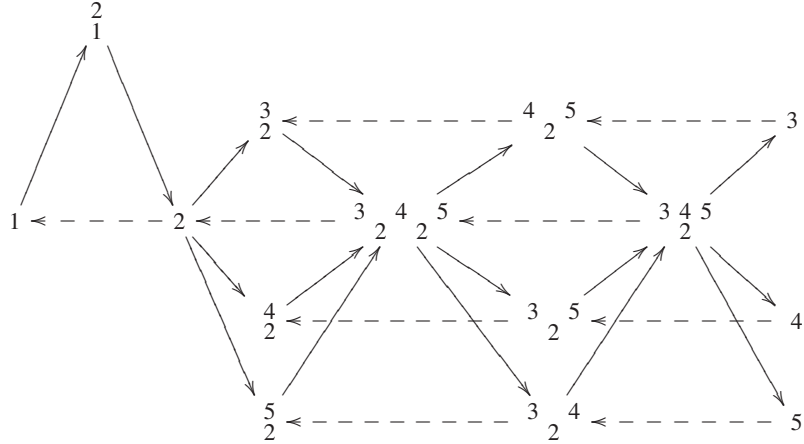
\xymatrix@C=1.5em@R=1.5em{
&
{\text{$\begin{smallmatrix}
2\\
1
\end{smallmatrix}$}}\ar[ddr]
&&&&&&\\
&&&
{\text{$\begin{smallmatrix}
3\\
2
\end{smallmatrix}$}}\ar[dr]
&&
{\text{$\begin{smallmatrix}
4&&5\\
&2&
\end{smallmatrix}$}}\ar[dr]\ar@{-->}[ll]
&&
{\text{$\begin{smallmatrix}
3
\end{smallmatrix}$}}\ar@{-->}[ll]
\\
{\text{$\begin{smallmatrix}
1
\end{smallmatrix}$}}\ar[uur]
&&
{\text{$\begin{smallmatrix}
2
\end{smallmatrix}$}}\ar[dr]\ar[ddr]\ar[ur]\ar@{-->}[ll]
&&
{\text{
$\begin{smallmatrix}
3&&4&&5 \\
&2&&2&
\end{smallmatrix}$}}\ar[dr]\ar[ddr]\ar[ur]\ar@{-->}[ll]
&&
{\text{
$\begin{smallmatrix}
3&4&5 \\
&2&
\end{smallmatrix}$}}\ar[dr]\ar[ddr]\ar[ur]\ar@{-->}[ll]
&
\\
&&&
{\text{
$\begin{smallmatrix}
4 \\
2
\end{smallmatrix}$}}\ar[ur]
&&
{\text{
$\begin{smallmatrix}
3&&5 \\
&2&
\end{smallmatrix}$}}\ar[ur]\ar@{-->}[ll]
&&
{\text{
$\begin{smallmatrix}
4
\end{smallmatrix}$}}\ar@{-->}[ll]
\\
&&&
{\text{
$\begin{smallmatrix}
5\\
2
\end{smallmatrix}$}}\ar[uur]
&&
{\text{
$\begin{smallmatrix}
3&&4 \\
&2&
\end{smallmatrix}$}}\ar[uur]\ar@{-->}[ll]
&&
{\text{
$\begin{smallmatrix}
5
\end{smallmatrix}$}}\ar@{-->}[ll]
}$$
Thus $A$ is representation-finite and $\ed(A)=0$ by Lemma \ref{ext-lem}(1). Let
$T:=
2\oplus \begin{matrix}2\\1\end{matrix}
\oplus \begin{matrix}3\\2\end{matrix}
\oplus \begin{matrix}4\\2\end{matrix}
\oplus \begin{matrix}5\\2\end{matrix}$.
By calculation, we obtain that $T$ is a tilting module with $\pd(_AT)=1$ . Let $P^{\bullet}$ be the projective resolution of $T$. Then $P^{\bullet}$ is a $2$-term tilting complex and $B:=\End_{\Db{A}}(P^{\bullet})$ is the path algebra given by the quiver $Q_B$
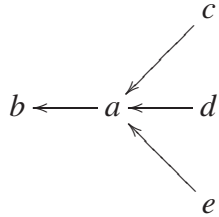
$$\xymatrix{
&&c\ar[dl]\\
b&a\ar[l]&d\ar[l]\\
&&e\ar[ul]
}$$
Since the underlying graph of $Q_B$ is Euclidean, we know that $B$ is a hereditary $k$-algebra of infinite representation type and $\ed(B)=1$ by Lemmas \ref{ext-lem}(1) and \ref{ext-lem}(3).

By the definition of $Q^{\bullet}$ in Lemma \ref{equi}, we can get $Q^{\bullet}$ is the $2$-term tilting complex over $B$ given by the direct sums of the following two complexes
$$\xymatrix{0\lra b\lra \text{$\begin{matrix}
a\\b
\end{matrix}$}}\lra 0,\quad
\xymatrix{0\lra \text{$
b\oplus
\begin{matrix}
c\\a\\b
\end{matrix}\oplus
\begin{matrix}
d\\a\\b
\end{matrix}\oplus
\begin{matrix}
e\\a\\b
\end{matrix}
$}}\lra 0\lra 0,$$
and $A\simeq \End_{\Db{B}}(Q^{\bullet})$ as algebras. Since $\gd(B)=1$, we have $\Ext^{2}_{B}(\T(Q^{\bullet}),\F(Q^{\bullet}))=0$. By Lemma \ref{sep}(1), $Q^{\bullet}$ is splitting and $P^{\bullet}$ is separating.

$(1)$ For the $2$-term silting complex $P^{\bullet}$, an easy calculation shows that $\F(P^{\bullet})=\add(_AS(1))$ and $\id(_AS(1))=2$. Since $\gd(B)=1$, we have $\pd(_BY)\le 1$ for each $Y\in \T(Q^{\bullet})$. By Lemma \ref{sep}(3), we know that $P^{\bullet}$ is separating. Due to $\ed(A)=0\neq 1=\ed(B)$, this example shows that the homological dimension restriction on $\F(P^{\bullet})$ cannot be removed.

$(2)$ For the $2$-term silting complex $Q^{\bullet}$, ($\T(Q^{\bullet})$, $\F(Q^{\bullet})$) in $B\modcat$ is not split, namely it is not separating. Although $\id(_BY)\le 1$ for every $Y\in \F(Q^{\bullet})$, $\ed(A)=0\neq 1=\ed(B)$. This implies that the separability condition in Theorem \ref{theo-silt} is necessary.

$(3)$ Note that $P^{\bullet}$ is $2$-term tilting complex and $|\ed(A)-\ed(B)|=1$. Thus this example also illustrates that the inequality in Theorem \ref{der-thm} can be an equality.
}\end{Bsp}

To illustrate Theorem \ref{theo-silt}, we give the following example. In addition, this example also implies that the inequality in Theorem \ref{der-thm} could be strict.
\begin{Bsp}{\rm
Let $k$ be a field, and $A$ be the path algebra given by the quiver
$$\xymatrix{
1\ar[r]^{\alpha}&2&3\ar[l]_{\beta}
}.$$
The Auslander-Reiten quiver is given by
$$
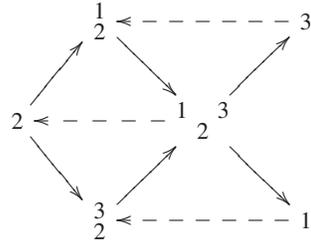
\xymatrix@C=1.5em@R=1.5em{
&
{\text{$\begin{smallmatrix}
1\\2
\end{smallmatrix}$}}\ar[dr]
&&
{\text{$\begin{smallmatrix}
3
\end{smallmatrix}$}}\ar@{-->}[ll]\\
{\text{$\begin{smallmatrix}
2
\end{smallmatrix}$}}\ar[ur]\ar[dr]
&&
{\text{$\begin{smallmatrix}
1&&3\\
&2&
\end{smallmatrix}$}}\ar[ur]\ar[dr]\ar@{-->}[ll]
&\\
&
{\text{$\begin{smallmatrix}
3\\2
\end{smallmatrix}$}}\ar[ur]
&&
{\text{$\begin{smallmatrix}
1
\end{smallmatrix}$}}\ar@{-->}[ll]
}$$
Let $P^{\bullet}$ be a complex given by the direct sums of the following two complexes{\small
$$0\lra 2\lra
\begin{matrix}
1\\2
\end{matrix}
\lra 0,\quad
0\lra 2\oplus
\begin{matrix}
3\\2
\end{matrix}
\lra 0\lra 0.$$}
By calculation, $P^{\bullet}$ is a $2$-term tilting complex over $A$, and $B:=\End_{\Db{A}}(P^{\bullet})$ is a path algebra given by the quiver
$$\xymatrix{
a&b\ar[l]&c\ar[l]
}.$$ The Auslander-Reiten quiver is given by
$$
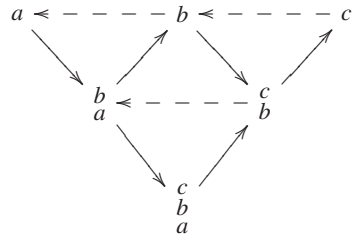
\xymatrix@C=1.5em@R=1.5em{
{\text{$\begin{smallmatrix}
a
\end{smallmatrix}$}}\ar[dr]
&&
{\text{$\begin{smallmatrix}
b
\end{smallmatrix}$}}\ar[dr]\ar@{-->}[ll]
&&
{\text{$\begin{smallmatrix}
c
\end{smallmatrix}$}}\ar@{-->}[ll]\\
&
{\text{$\begin{smallmatrix}
b\\a
\end{smallmatrix}$}}\ar[ur]\ar[dr]
&&
{\text{$\begin{smallmatrix}
c\\b
\end{smallmatrix}$}}\ar[ur]\ar@{-->}[ll]
&\\
&&
{\text{$\begin{smallmatrix}
c\\b\\a
\end{smallmatrix}$}}\ar[ur]
&&
}$$
By direct computation, $Q^{\bullet}$ (defined in Lemma \ref{equi}) is the complex given by the direct sums of the following complexes
{\small
$$0\lra a\lra
\begin{matrix}
b\\a
\end{matrix}\lra 0,\quad
0\lra 0\lra \begin{matrix}
b\\a
\end{matrix}\oplus
\begin{matrix}
c\\ b \\ a
\end{matrix}\lra 0.
$$}
Further, the induced torsion pair $(\T(P^{\bullet}),\F(P^{\bullet}))$ in $A\modcat$ is given by
$$
\T(P^{\bullet}) = \add(
\begin{smallmatrix}
1
\end{smallmatrix}),\quad
\F(P^{\bullet}) = \add(
\begin{smallmatrix}
2
\end{smallmatrix}\oplus
\begin{smallmatrix}
1\\2
\end{smallmatrix}\oplus
\begin{smallmatrix}
3\\2
\end{smallmatrix}\oplus
\begin{smallmatrix}
1&&3\\
&2&
\end{smallmatrix}\oplus
\begin{smallmatrix}
3
\end{smallmatrix}).$$
The induced torsion pair $(\T(Q^{\bullet}),\F(Q^{\bullet}))$ in $B\modcat$ is given by
$$\T(Q^{\bullet})=\add(
\begin{matrix}
b\\a
\end{matrix}\oplus
\begin{matrix}
b
\end{matrix}\oplus
\begin{matrix}
c\\ b \\ a
\end{matrix}\oplus
\begin{matrix}
c\\b
\end{matrix}\oplus
\begin{matrix}
c
\end{matrix}),\quad
\F(Q^{\bullet})=\add(\begin{matrix}
a
\end{matrix}).
$$
Clearly, $(\T(P^{\bullet}),\F(P^{\bullet}))$ and $(\T(Q^{\bullet}),\F(Q^{\bullet}))$ are split. Then the homological conditions in Theorem \ref{theo-silt} is satisfied. Thus $\ed(A)=\ed(B)$ and $|\ed(A)-\ed(B)|=0<1$. This implies the inequality in Theorem \ref{der-thm} can be strict.
}\end{Bsp}

\section{Stable equivalences}
In this section, we shall prove that the extension dimension is invariant under stable equivalence. We first recall some basic results about the stable equivalence of Artin algebras , as detailed in reference \cite{ars97,c21,Guo05}.

Let $A$ be an Artin algebra over a fixed commutative Artin ring $k$. Recall that an $A$-module $X$ is called a \emph{generator} if $A\in\add(X)$, \emph{cogenerator} if $\DD(A_A)\in\add(X)$, and \emph{generator-cogenerator} if it is both a generator and cogenerator for $A\modcat$. We denoted by $\stmodcat{A}$ the stable module category of $A$ modulo projective modules. The objects are the same as the objects of $A\modcat$, and for two modules $X, Y$ in $\stmodcat{A}$, their homomorphism set is $\underline{\Hom}_A(X,Y):=\Hom_A(X,Y)/\mathscr{P}(X,Y)$, where $\mathscr{P}(X,Y)$ is the subgroup of $\Hom_A(X,Y)$ consisting of the homomorphisms factorizing through a projective $A$-module. This category is usually called the \emph{stable module category} of $A$. Dually, We denoted by  $A\mbox{-}\overline{\mbox{mod}}$ the stable module category of $A$ modulo injective modules. Let $\tau_A$ be the Auslander-Reiten translation $\DTr$. Then $\tau_A:\stmodcat{A}\ra A\mbox{-}\overline{\mbox{mod}}$ be an equivalence as additive categories (see \cite[Chapler IV.1]{ars97}). Two algebras $A$ and $B$ are said to be \emph{stably equivalent} if the two stable categories $A\stmodcat$ and $B\stmodcat$ are equivalent as additive categories.

Next, suppose that $F: \stmodcat{A} \ra \stmodcat{B}$ is an stable equivalence. Then the following functor
$$F':=\tau_B\circ F\circ \tau_A^{-1}:A\mbox{-}\overline{\mbox{mod}}\lra B\mbox{-}\overline{\mbox{mod}}$$
is equivalent as additive categories. Moreover, there are one-to-one correspondences
$$F: A\modcat_{\mathscr{P}}\lra B\modcat_{\mathscr{P}}\quad \mbox{and}\quad
F': A\modcat_{\mathscr{I}}\lra  B\modcat_{\mathscr{I}},$$
where $A\modcat_{\mathscr{P}}$ (respectively, $A\modcat_{\mathscr{I}}$) stands for the full subcategory of $A\modcat$ consisting of modules without nonzero projective (respectively, injective) summands.
We also use $F$ (respectively, $F'$) to denote the induce map $A\modcat\ra B\modcat$ which takes projectives (respectively, injectives) to zero.

Recall that a simple $A$-module $S$ is called a \emph{node} of $A$ if it is  neither projective nor injective, and the middle term of the almost split sequence starting at $S$ is projective; a node $S$ in $A\modcat$ is said to be an $F$-\emph{exceptional node} if $F(S)\not\simeq F'(S)$. Let $\mathfrak{n}_{F}(A)$ be the set of isomorphism classes of $F$-exceptional nodes of $A$. Since $\mathfrak{n}_{F}(A)$ is a subset of all simple modules, $\mathfrak{n}_{F}(A)$ is a finite set. Note that $X$ is indecomposable, non-projective, non-injective, and not a node, then $F(X)\simeq F'(X)$ (see \cite[Lemma 3.4]{ar78} or \cite[Chapter X.1.7, p. 340]{ars97}). Then $\mathfrak{n}_{F}(A)$ and the set of isomorphism classes
of indecomposable, non-projective, non-injective $A$-modules $U$ such that $F(U)\not\simeq F'(U)$, are coincide.

Let $F^{-1}:\stmodcat{B}\ra \stmodcat{A}$ be a quasi-inverse of $F$. Then we use $\mathfrak{n}_{F^{-1}}(B)$ to denote the the set of isomorphism classes of $F^{-1}$-exceptional nodes of $B$.

In the following, let $$\bigtriangleup_A
:=\mathfrak{n}_{F}(A)
\dot{\cup}(\mathscr{P}_A\setminus\mathscr{I}_A)
\quad\mbox{and}\quad \bigtriangledown_A
:=\mathfrak{n}_{F}(A)\dot{\cup}(\mathscr{I}_A\setminus\mathscr{P}_A),$$ where $\dot{\cup}$ stands for the disjoint union of sets; $\mathscr{P}_A$ and $\mathscr{I}_A$ stand for the set of isomorphism classes of indecomposable projective and injective $A$-modules, respectively. By $\bigtriangleup_A^c$ we mean the class of indecomposable, non-injective $A$-modules which do not belong to $\bigtriangleup_A$.

\begin{Rem}\label{rem-tri}
Each module $X\in A\modcat_{\mathscr{I}}$ admits a unique decomposition (up to isomorphism)
$$X\simeq X^{c}\oplus X^{\triangle}$$
with $X^{c}\in \bigtriangleup_A^c$ and $X^{\triangle}\in \bigtriangleup_A$.
\end{Rem}
Let $\mathcal{GCN}_{F}(A)$ be
the class of basic $A$-modules $X$ which are generator-cogenerators with  $\mathfrak{n}_{F}(A)\subseteq \add(X)$, that is,
$$\mathcal{GCN}_{F}(A)=\{M'\oplus M_0\mid M'\in A\modcat \mbox{ and }M' \mbox{ is basic} \},$$
where $M_0$ is the unique (up to isomorphism) basic module with $$\add(M_0)=\add(A)\cup\add(\DD(A_A))\cup\add(\mathfrak{n}_{F}(A)).$$
We say a module is called \emph{basic} if it is a direct sum of non-isomorphic indecomposable modules.

We define the following correspondences:
$$\Phi:A\modcat\lra B\modcat,\quad X\mapsto F(X)\oplus\bigoplus_{Q\in \mathscr{P}_B}Q,$$
$$\Psi:B\modcat\lra A\modcat,\quad Y\mapsto F^{-1}(Y)\oplus\bigoplus_{P\in \mathscr{P}_A}P.$$

\begin{Lem}\label{one-to-one}{\rm (\cite[Lemma 4.10]{c21})}
$(1)$ There exist one-to-one correspondences
$$F:\bigtriangledown_A\lra \bigtriangledown_B
,\quad
F':\bigtriangleup_A\lra \bigtriangleup_B
\quad\mbox{and}\quad
F':\bigtriangleup_A^c\lra \bigtriangleup_B^c.$$

$(2)$ The correspondences $\Phi$ and $\Psi$ restrict to one-to-one correspondences between $\mathcal{GCN}_{F}(A)$ and $\mathcal{GCN}_{F^{-1}}(B)$. Moreover, if $X\in \mathcal{GCN}_{F}(A)$, then $\Phi(X)\simeq F'(X)\oplus \bigoplus_{J\in\mathscr{I}_B}J$.
\end{Lem}

Recall that an exact sequence $0\ra X\lraf{f} Y\ra Z\lraf{g} 0$ in $A\modcat$ is called \emph{minimal} if it has no a split exact sequence as a direct summand, that is, there does not exist isomorphisms $u$, $v$, $w$ such that the following diagram
$$\xymatrix{
0\ar[r]
&X\ar[r]^{f}\ar[d]^{u}
&Y\ar[r]^g\ar[d]^{v}
&Z\ar[r]\ar[d]^{w}
&0\\
0\ar[r]&X_1\oplus
X_2\ar[r]^{\text{
$\left(
\begin{smallmatrix}
f_1&0\\
0&f_2
\end{smallmatrix}\right)
$}}
&Y_1\oplus Y_2\ar[r]^{\text{
$\left(
\begin{smallmatrix}
g_1&0\\
0&g_2
\end{smallmatrix}\right)
$}}&Z_1\oplus Z_2\ar[r]&0
}$$
is row exact and commute, where $0\ra X_2\lraf{f_2} Y_2\lraf{g_2}Z_2\ra 0$ is split. The next lemma shows that the stable functor has certain ``exactness" property.
\begin{Lem}\label{exact}
Let $Z$ be an $A$-module without nonzero projective summands, and let $$0\lra X\oplus X'\lra Y\oplus P\lraf{g}Z\lra 0$$ be a minimal short exact sequence in $A\modcat$ such that $X\in \add(\bigtriangleup_A^c)$, $X'\in \add(\bigtriangleup_A)$, $Y\in A\modcat_{\mathscr{P}}$ and $P\in\add(_AA)$.
Then there exists a minimal short exact sequence $$0\lra F(X)\oplus F'(X')\lra F(Y)\oplus Q\lraf{g'}F(Z)\lra 0$$ in $B\modcat$
such that $Q\in \add(_BB)$ and $g'=F(g)$ in $\stmodcat{B}$.
\end{Lem}
{\it Proof.} Note that if $_AZ$ is indecomposable and non-projective, then the statement in Lemma \ref{exact} has been proved in \cite[Lemma 4.13]{c21}, which is valid also for decomposable module having no nonzero projective summands by checking the argument there. $\square$

\begin{Lem}\label{mwrd}
Let $M\in \mathcal{GCN}_{F}(A)$. Then $M$-$\wrd(A)=\Phi(M)$-$\wrd(B)$.
\end{Lem}
{\it Proof.} Suppose $M$-$\wrd(A)=n$. We first prove that $\Phi(M)$-$\wrd(B)\le n$. For $Y\in B\modcat$, we can write $Y\simeq Y^{\mathscr{P}}\oplus Q'$ with $Y^{\mathscr{P}}\in B\modcat_{\mathscr{P}}$ and $Q'\in \pmodcat{B}$. Since $F: A\modcat_{\mathscr{P}}\ra B\modcat_{\mathscr{P}}$ is an one-to-one correspondence, there exists $X\in A\modcat_{\mathscr{P}}$ such that $F(X)\simeq Y^{\mathscr{P}}$ as $B$-modules. By Definition \ref{wrd-def}, there exists an exact sequence
$$0\lra M_n\lraf{f_n} M_{n-1} \lraf{f_{n-1}} \cdots \lraf{f_1} M_0 \lraf{f_0} X\oplus X'\lra 0$$
in $A\modcat$ for some $A$-module $X'$ such that $M_i\in\add(_AM)$ for $0\le i\le n-1$. Then we have short exact sequences
\begin{align}\label{kmk}
0\lra K_{i+1}\lra M_i\lra K_i\lra 0
\end{align}
where $K_0:=X\oplus X'$ and $K_{i+1}$ is the kernel of $f_i$ for each $0\le i\le n-1$.

For $0\le i\le n-1$, we can decompose the short exact sequence (\ref{kmk}) as the direct sums of the following two short exact sequences
\begin{align}\label{uwu}
0\lra U_{i+1}\lra W_i\lra U_i\lra 0,
\end{align}
\begin{align}\label{vzv}
0\lra V_{i+1}\lra Z_i\lra V_i\lra 0
\end{align}
in $A\modcat$, namely there are isomorphisms $\lambda_i$, $\mu_i$, $\nu_i$ with the following commutative diagram
\begin{align}\label{2row}
\xymatrix@C=0.5em@R=0.5em{
0\ar[rr]
&&K_{i+1}\ar[rr]\ar[dd]^{\lambda_i}
&&M_i\ar[rr]\ar[dd]^{\mu_i}
&&K_i\ar[rr]\ar[dd]^{\nu_i}
&&0\\
&&&&&&&&\\
0\ar[rr]
&&U_{i+1}\oplus
V_{i+1}\ar[rr]
&&W_i\oplus Z_i\ar[rr]
&&U_i\oplus V_i\ar[rr]
&&0,
}\end{align}
in $A\modcat$ such that (\ref{uwu}) is minimal and (\ref{vzv}) is split. Since the sequence (\ref{uwu}) is minimal, we have $U_{i+1}\in A\modcat_{\mathscr{I}}$.
By Remark \ref{rem-tri}, we can write $U_{i+1}\simeq U_{i+1}^{c}\oplus U_{i+1}^{\triangle}$
with $U_{i+1}^{c}\in \bigtriangleup_A^c$ and $U_{i+1}^{\triangle}\in \bigtriangleup_A$.
We also write $W_i:=W^{^{\mathscr{P}}}_i\oplus P_i$ with $W^{\mathscr{P}}_i\in A\modcat_{\mathscr{P}}$ and $P_i\in \pmodcat{A}$.
Then we can write the following short exact sequence
\begin{align}\label{uuu}
0\lra U_{i+1}^{c}\oplus U_{i+1}^{\triangle}\lra W^{\mathscr{P}}_i\oplus P_i\lra U_i\lra 0
\end{align}
for the sequence (\ref{uwu}). Since the sequence (\ref{uwu}) is minimal, we know $U_{i} \in A\modcat_{\mathscr{P}}$. By Lemma \ref{exact}, we have the following minimal exact sequence
$$ 0\lra F(U_{i+1}^{c})\oplus F'(U_{i+1}^{\triangle})\lra F(W^{\mathscr{P}}_i)\oplus Q_i\lra F(U_i)\lra 0$$
in $B\modcat$ such that $Q_i\in \add(_BB)$ for each $0\le i\le n-1$, that is, we have the following exact sequences
\begin{equation*}
\begin{cases}
\xymatrix@C=1em@R=0.1em{
0\ar[r]& F(U_{1}^{c})\oplus F'(U_{1}^{\triangle})\ar[r]& F(W^{\mathscr{P}}_0)\oplus Q_0\ar[r]& F(U_0)\ar[r] &0\\
0\ar[r]& F(U_{2}^{c})\oplus F'(U_{2}^{\triangle})\ar[r]& F(W^{\mathscr{P}}_1)\oplus Q_1\ar[r]& F(U_{1}^{c})\oplus F(U_{1}^{\triangle})\ar[r] &0\\
0\ar[r]& F(U_{3}^{c})\oplus F'(U_{3}^{\triangle})\ar[r]& F(W^{\mathscr{P}}_2)\oplus Q_2\ar[r]& F(U_{2}^{c})\oplus F(U_{2}^{\triangle})\ar[r] &0\\
&&\vdots&&\\
0\ar[r]& F(U_{n-1}^{c})\oplus F'(U_{n-1}^{\triangle})\ar[r]& F(W^{\mathscr{P}}_{n-2})\oplus Q_{n-2}\ar[r]& F(U_{n-2}^{c})\oplus F(U_{n-2}^{\triangle})\ar[r] &0\\
0\ar[r]& F(U_{n}^{c})\oplus F'(U_{n}^{\triangle})\ar[r]& F(W^{\mathscr{P}}_{n-1})\oplus Q_{n-1}\ar[r]& F(U_{n-1}^{c})\oplus F(U_{n-1}^{\triangle})\ar[r] &0.}
\end{cases}
\end{equation*}
Furthermore, we have short exact sequences
{\footnotesize
\begin{equation}\label{short-st}
\begin{cases}
\xymatrix@C=1em@R=0.1em{
0\ar[r]
& F(U_{1}^{c})\oplus F'(U_{1}^{\triangle})\oplus F(U_{1}^{\triangle})\ar[r]
& F(W^{\mathscr{P}}_0)\oplus Q_0\oplus F(U_{1}^{\triangle})\oplus F(V_0)\ar[r]
& F(U_0)\oplus F(V_0)\ar[r]
&0\\
0\ar[r]
& F(U_{2}^{c})\oplus F'(U_{2}^{\triangle})\oplus F(U_{2}^{\triangle})\ar[r]
& F(W^{\mathscr{P}}_1)\oplus Q_1\oplus F'(U_{1}^{\triangle})\oplus F(U_{2}^{\triangle})\ar[r]
& F(U_{1}^{c})\oplus F(U_{1}^{\triangle})\oplus F'(U_{1}^{\triangle})\ar[r]
&0\\
0\ar[r]
& F(U_{3}^{c})\oplus F'(U_{3}^{\triangle})\oplus F(U_{3}^{\triangle})\ar[r]
& F(W^{\mathscr{P}}_2)\oplus Q_2\oplus F'(U_{2}^{\triangle})\oplus F(U_{3}^{\triangle})\ar[r]
& F(U_{2}^{c})\oplus F(U_{2}^{\triangle})\oplus F'(U_{2}^{\triangle})\ar[r]
&0\\
&&\vdots&&\\
0\ar[r]
& F(U_{n-1}^{c})\oplus F'(U_{n-1}^{\triangle})\oplus F(U_{n-1}^{\triangle})\ar[r]
& F(W^{\mathscr{P}}_{n-2})\oplus Q_{n-2}\oplus F'(U_{n-2}^{\triangle})\oplus F(U_{n-1}^{\triangle})\ar[r]
& F(U_{n-2}^{c})\oplus F(U_{n-2}^{\triangle})\oplus F'(U_{n-2}^{\triangle})\ar[r]
&0\\
0\ar[r]
& F(U_{n}^{c})\oplus F'(U_{n}^{\triangle})\ar[r]
& F(W^{\mathscr{P}}_{n-1})\oplus Q_{n-1}\oplus F'(U_{n-1}^{\triangle})\ar[r]
& F(U_{n-1}^{c})\oplus F(U_{n-1}^{\triangle})\oplus F'(U_{n-1}^{\triangle})\ar[r]
&0.}
\end{cases}
\end{equation}}
Combine the short exact sequences (\ref{short-st}) ,we get a long exact sequence in $B\modcat$:
\begin{align}\label{long-st}
\xymatrix{
&0\lra F(U_{n}^{c})\oplus F'(U_{n}^{\triangle}) \ar[r]
& F(W^{\mathscr{P}}_{n-1})\oplus Q_{n-1}\oplus F'(U_{n-1}^{\triangle}) \ar `r[d] `[l] `[lld] `[dl] [dl]\\
& F(W^{\mathscr{P}}_{n-2})\oplus Q_{n-2}\oplus F'(U_{n-2}^{\triangle})\oplus F(U_{n-1}^{\triangle}) \ar[r]
& \cdots \lra F(W^{\mathscr{P}}_1)\oplus Q_1\oplus F'(U_{1}^{\triangle})\oplus F(U_{2}^{\triangle})
\ar `r[d] `[l] `[lld] `[dl] [dl]\\
& F(W^{\mathscr{P}}_0)\oplus Q_0\oplus F(U_{1}^{\triangle})\oplus F(V_0) \ar[r]
&  F(U_0)\oplus F(V_0)  \lra 0.
}\end{align}
By the commutative diagram (\ref{2row}), we get $M_n=K_n\simeq U_n\oplus V_n$ and $M_i\simeq W_i\oplus Z_i$ for $0\le i\le n-1$.
Also $U_n\simeq U_n^{c}\oplus U_n^{\triangle}$ and $W_i\simeq W_i^{\mathscr{P}}\oplus Q_i$ for $0\le i\le n-1$.
Thus $U_n^{c},W_i^{\mathscr{P}},Z_i\in \add(_AM)$ for $0\le i\le n-1$. Since $M\in \mathcal{GCN}_{F}(A)$, we have $\add(_AA)\cup \add(_A\DD(A_A))\cup\add(_A\mathfrak{n}_{F}(A))\subseteq\add(_AM)$.
In particular, $\add(\bigtriangleup_A)\subseteq \add(_AM)$.
Thus $U^{\triangle}_i\in \add(_AM)$ for $1\le i\le n$.
Since the exact sequence $0\ra V_1\ra Z_0\ra V_0\ra 0$ (see (\ref{vzv})) is split, we have $V_0\in \add(_AZ_0)\subseteq \add(_AM)$.
By Lemma \ref{one-to-one}, $$\Phi(M)=F(M)\oplus\bigoplus_{Q\in \mathscr{P}_B}Q\simeq F'(M)\oplus \bigoplus_{J\in\mathscr{I}_B}J.$$
Thus $F(U_n^{c}),F(W_i^{\mathscr{P}}),F(U^{\triangle}_{i+1}),F'(U^{\triangle}_{i+1}),F(V_0)\in \add(_B\Phi(M))$ for $0\le i\le n-1$.
Also $Q_i\in \add(_B\Phi(M))$.
By the long exact sequence (\ref{long-st}) and Definition \ref{wrd-def}, we get $$\Phi(M)\mbox{-}\wrd(_B(F(U_0)\oplus F(V_0)))\le n.$$
Due to the commutative diagram (\ref{2row}), we have $K_0\simeq U_0\oplus V_0$ as $A$-modules.
It follows from $K_0=X\oplus X'$ that $F(X)\oplus F(X')\simeq F(U_0)\oplus F(V_0)$ as $B$-modules.
Thus $$\Phi(M)\mbox{-}\wrd(_BF(X))\le \Phi(M)\mbox{-}\wrd(_B(F(U_0)\oplus F(V_0)))\le n.$$
Note that $Y\simeq Y^{\mathscr{P}}\oplus Q'$ as $B$-modules with $Y^{\mathscr{P}}\simeq F(X)$ and $Q'\in \add(_BB)\subseteq \add(_B\Phi(M))$.
Thus $$\Phi(M)\mbox{-}\wrd(_BY)=\Phi(M)\mbox{-}\wrd(_BF(X))\le n= M\mbox{-}\wrd(A).$$
By Definition \ref{wrd-def}, we have $$\Phi(M)\mbox{-}\wrd(B)\le M\mbox{-}\wrd(A).$$ Similarly, we get $$\Phi(M)\mbox{-}\wrd(B)\ge M\mbox{-}\wrd(A).$$ Thus $\Phi(M)\mbox{-}\wrd(B)= M\mbox{-}\wrd(A).$
$\square$

\begin{Theo}\label{st-thm}
Let $A$ and $B$ be stably equivalent Artin algebras. Then $\ed(A)=\ed(B)$.
\end{Theo}
{\it Proof.} Let $M_0$ be the unique (up to isomorphism) basic module with $$\add(_AM_0)=\add(A)\cup\add(\DD(A_A))\cup\add(\mathfrak{n}_{F}(A)).$$
Then $\mathcal{GCN}_{F}(A)=\{M'\oplus M_0\mid M'\in A\modcat \mbox{ and }M' \mbox{ is basic} \}.$
By Lemma \ref{weak-prop}(3), we have
$$\min\{M \mbox{-}\wrd(A)\mid M\in A\modcat\}
=\min\{M \mbox{-}\wrd(A)\mid M\in\mathcal{GCN}_{F}(A)\}.$$
Similarly, we have
$$\min\{N \mbox{-}\wrd(B)\mid N\in \mathcal{GCN}_{F^{-1}}(B)\}
=\min\{N \mbox{-}\wrd(B)\mid N\in B\modcat\}.$$
By Definition \ref{wrd-def}, we get the following equalities.
\begin{align*}
\wrd(A)
&=\min\{M \mbox{-}\wrd(A)\mid M\in A\modcat\}\\
&=\min\{M \mbox{-}\wrd(A)\mid M\in \mathcal{GCN}_{F}(A)\}\\
&=\min\{N \mbox{-}\wrd(B)\mid N\in \mathcal{GCN}_{F^{-1}}(B)\}
\quad (\mbox{ by Lemmas \ref{one-to-one}(2) and \ref{mwrd} })\\
&=\min\{N \mbox{-}\wrd(B)\mid N\in B\modcat\}\\
&=\wrd(B).
\end{align*}
By Lemma \ref{wrd=ed}, we have $\ed(A)=\ed(B)$.
$\square$

\medskip
Now, we deduce some consequences of Theorem \ref{st-thm}.

Recall that a derived equivalence $F$ between finite dimensional algebras $A$ and $B$ with a quasi-inverse $G$ is called \emph{almost $\nu$-stable} \cite{hx10} if the associated radical tilting complexes $\cpx{T}$ over $A$ to $F$ and $\cpx{\bar{T}}$ over $B$ to $G$ are of the form
$$ \cpx{T}: 0\lra T^{-n}\lra \cdots T^{-1}\lra T^0\lra 0 \; \mbox{ and } \,
\cpx{\bar{T}}: 0\lra \bar{T}^{0}\lra \bar{T}^1\lra \cdots \lra \bar{T}^{n}\lra 0,$$
respectively, such that $\add(\bigoplus_{i=1}^n T^{-i})=\add(\nu_A(\bigoplus_{i=1}^n T^{-i}))$ and $\add(\bigoplus_{i=1}^n \bar{T}^{i}) = \add(\nu_B(\bigoplus_{i=1}^n \bar{T}^{i}))$, where $\nu$ is the Nakayama functor. By \cite[Theorem 1.1(2)]{hx10}, almost $\nu$-stable derived equivalences induce special stable equivalences, namely stable equivalences of Morita type. Thus we have the following consequence of Theorem \ref{st-thm}.

\begin{Koro}\label{nu-s}
Let $A$ and $B$ be almost $\nu$-stable derived equivalent finite dimensional algebras. Then $\ed(A)=\ed(B)$.
\end{Koro}

Recall that given a finite dimensional algebra $A$ over a filed $k$, $A\ltimes D(A)$, the trivial extension of $A$ by $D(A)$ is the $k$-algebra whose underlying $k$-space is $A\oplus D(A)$, with multiplication given by

$$ (a,f)(b, g)= (ab,fb+ag)$$ for $a,b\in A$, and $f,g\in D(A)$, where $D:=\Hom_{k}(-,k)$. It is known that $A\ltimes D(A)$ is always symmetric, and therefore it is selfinjective.

\begin{Koro}
Let $A$ and $B$ be derived equivalent finite dimensional algebras. Then $$\ed(A\ltimes D(A))=\ed(B\ltimes D(B)).$$
\end{Koro}
{\it Proof.} By a result of Rickard (see \cite[Theorem 3.1]{r89}), which says that any derived equivalence between two algebras induces a derived equivalence between their trivial extension algebras, we obtain that $A\ltimes D(A)$ and $B\ltimes D(B)$ are derived equivalent.  It follows from \cite[Proposition 3.8]{hx10} that every derived equivalence between two selfinjective algebras induces an almost $\nu$-stable derived equivalence. Thus we have $\ed(A\ltimes D(A))=\ed(B\ltimes D(B))$ by Corollary \ref{nu-s}.
$\square$

\medskip
In general, it is rather hard to compute the precise value of the extension dimension of a given algebra. However, we display an example to show how Theorem \ref{st-thm} can be applied to compute the extension dimensions of certain algebras. The example shows also that the method of computing extension dimensions by stable equivalences seems to be useful.
\begin{Bsp}{\rm
Let $k$ be a fixed field, $A=kQ_A/I$, where $Q_A$ is the quiver
$$\xymatrix@R=10pt{
&&1\ar[d]^{\beta}\ar@(ur,ul)_{\gamma}&&&&\\
2\ar[r]^{\alpha_2}
&3\ar[r]^{\alpha_3}
&4\ar[r]^{\alpha_4}
&5\ar[r]^{\alpha_5}
&6\ar[r]^{\alpha_6}
&\cdots \ar[r]^{\alpha_{n-1}}
&n
}$$
and $I$ is generated by $\{\gamma^2, \beta\gamma\}$ with $n\ge 6$.
The indecomposable projective left $A$-modules are
$$
\xymatrix@R=.23cm@C=.01cm{
P(1)=&&1\ar@{-}[dl]\ar@{-}[dr]&\\
&1&&4\ar@{-}[d]\\
&&&5\ar@{-}[d]\\
&&&6\ar@{-}[d]\\
&&&\vdots\ar@{-}[d]\\
&&&n}\quad
\xymatrix@R=.23cm@C=.01cm{
P(2)=&2\ar@{-}[d]\\
&3\ar@{-}[d]\\
&4\ar@{-}[d]\\
&5\ar@{-}[d]\\
&\vdots\ar@{-}[d]\\
&n}\quad
\xymatrix@R=.23cm@C=.01cm{
P(3)=&3\ar@{-}[d]\\
&4\ar@{-}[d]\\
&5\ar@{-}[d]\\
&6\ar@{-}[d]\\
&\vdots\ar@{-}[d]\\
&n}\quad
\cdots\quad
\xymatrix@R=.23cm@C=.01cm{
P(n-2)=&n-2\ar@{-}[d]\\
&n-1\ar@{-}[d]\\
&n}
\quad
\xymatrix@R=.23cm@C=.01cm{
P(n-1)=&n-1\ar@{-}[d]\\
&n}
\quad
\xymatrix@R=.23cm@C=.01cm{
P(n)=&n
}$$
Then we have $\pd(S(1))=\infty$, $\pd(S(i))=1$ for $2\le i\le n-1$ and $\pd(S(n))=0$. Clearly, $\gd(A)=\infty$.
By calculation, we get that $\ell\ell(A)=n-1$ and $\ell\ell^{\infty}(A)=2$, where $\ell\ell(A)$ and $\ell\ell^{\infty}(A)$ stand for the Loewy length and the infinite-layer length of $A$ {\rm(\cite{hlh09,hlh13})}, respectively. By Lemma \ref{ext-lem}, we have $$\ed(A)\le \min\{\gd(A),\ell\ell(A)-1,\ell\ell^{\infty}(A)+1\}=3.$$

By the previous works in Lemma \ref{ext-lem}, one just get the upper bound for the extension dimension of $A$. Next, we shall compute the extension dimension of $A$ by our results. By \cite[Lemma 1]{MV1980}, we know that $S(1)$ is a unique node of $A$. It follows from \cite[Theorem 2.10]{MV1980} that $A$ is stably equivalent to the path algebra $B$ given by the following quiver $Q_B$:
$$\xymatrix@R=10pt{
&&1'&&&&\\
&&1\ar[d]^{\beta}\ar[u]_{\delta}&&&&\\
2\ar[r]^{\alpha_2}
&3\ar[r]^{\alpha_3}
&4\ar[r]^{\alpha_4}
&5\ar[r]^{\alpha_5}
&6\ar[r]^{\alpha_6}
&\cdots \ar[r]^{\alpha_{n-1}}
&n.}$$
Since the underlying graph of $Q_B$ is not Dynkin, $B$ is representation-infinite and $\ed(B)=1$ by Lemmas \ref{ext-lem}(1) and \ref{ext-lem}(3). By Theorem \ref{st-thm}, we have $\ed(A)=\ed(B)=1$.
}\end{Bsp}

\bigskip
\noindent{\bf Acknowledgements.}
Jinbi Zhang would like to thank Prof. Changchang Xi from Capital Normal University and Prof. Jiping Zhang from Peking University for help and encouragement. Junling Zheng was supported by the National Nature Since Foundation of China(Grant No. 12001508).

\end{document}